\documentclass[11pt]{amsart}
\usepackage{amsxtra}
\usepackage{amssymb}
\theoremstyle{plain}
\newcommand{\cleqn}{\setcounter{equation}{0}}
\newcommand{\clth}{\setcounter{theorem}{0}}
\newcommand {\sectionnew}[1]{\section{#1}\cleqn\clth}
\newcommand{\nn}{\hfill\nonumber}
\newtheorem{theorem}{Theorem}[section]
\newtheorem{lemma}[theorem]{Lemma}
\newtheorem{definition-theorem}[theorem]{Definition-Theorem}
\newtheorem{proposition}[theorem]{Proposition}
\newtheorem{corollary}[theorem]{Corollary}
\newtheorem{definition}[theorem]{Definition}
\newtheorem{example}[theorem]{Example}
\newtheorem{remark}[theorem]{Remark}
\newtheorem{conjecture}[theorem]{Conjecture}
\newtheorem{notation}[theorem]{Notation}
\newcommand \bth[1] { \begin{theorem}\label{t#1} }
\newcommand \ble[1] { \begin{lemma}\label{l#1} }

\newcommand \bpr[1] { \begin{proposition}\label{p#1} }
\newcommand \bco[1] { \begin{corollary}\label{c#1} }
\newcommand \bde[1] { \begin{definition}\label{d#1}\rm }
\newcommand \bex[1] { \begin{example}\label{e#1}\rm }
\newcommand \bre[1] { \begin{remark}\label{r#1}\rm }
\newcommand \bcj[1] { \begin{conjecture}\label{j#1}\rm }

\newcommand \bnota[1] { \begin{notation}\label{n#1}\rm }
\renewcommand {\eth} { \end{theorem} }
\newcommand {\ele} { \end{lemma} }

\newcommand {\epr} { \end{proposition} }
\newcommand {\eco} { \end{corollary} }
\newcommand {\ede} { \end{definition} }
\newcommand {\eex} { \end{example} }
\newcommand {\ere} { \end{remark} }
\newcommand {\ecj} { \end{conjecture} }

\newcommand {\enota} { \end{notation} }
\newcommand \thref[1]{Theorem \ref{t#1}}
\newcommand \leref[1]{Lemma \ref{l#1}}
\newcommand \prref[1]{Proposition \ref{p#1}}

\newcommand \lb[1]{\label{#1}}
\def \Cset {{\mathbb C}}
\def \KK {{\mathbb K}}
\def \Zset {{\mathbb Z}}
\def \Nset {{\mathbb N}}
\def \Qset {{\mathbb Q}}

\def \Tset {{\mathbb T}}

\def \B  {{\mathcal{B}}}

\def \QQ {{\mathcal{Q}}}
\def \PP {{\mathcal{P}}}

\def \VV {{\mathcal{V}}}
\def \UU {{\mathcal{U}}}
\def \RR {{\mathcal{R}}}
\def \SS {{\mathcal{S}}}
\def \LL {{\mathcal{L}}}
\def \TT {{\mathcal{T}}} 
\def \ZZ {{\mathcal{Z}}}
\def \pb {{\bf{p}}}
\def \rb {{\bf{r}}}
\def \De {\Delta}   % Greek letters
\def \de {\delta}
\def \al {\alpha}
\def \be {\beta}
\def \la {\lambda}

\def \om {\omega}
\def \Om {\Omega}
\def \ga {\gamma}
\def \de {\delta}
\def \Ga {\Gamma}

\def \sig {\sigma}

\def \del {\partial}
\def \ep {\epsilon}
\def \sig{\sigma}
\def \mt  {\mapsto}
\def \hra {\hookrightarrow}
\def \lha {\leftharpoonup}
\def \rha {\rightharpoonup}

\def \rcor {\rangle}
\def \lcor {\langle}
\def \o  {\otimes}
\def \ol {\overline}

\def \wh {\widehat}

\def \id { {\mathrm{id}} }

\def \rank { {\mathrm{rank}} }

\def \g  {\mathfrak{g}}   % Lie algebra letters

\def \n  {\mathfrak{n}}

\DeclareMathOperator \Span { {\mathrm{Span}} }

\DeclareMathOperator \charr { {\mathrm{char}} }

\DeclareMathOperator \Ker { {\mathrm{Ker}} }

\DeclareMathOperator \Ext { {\mathrm{Ext}} }
\DeclareMathOperator \GKdim {{\mathrm{GK \, dim}}}
\DeclareMathOperator \Wt  { {\mathrm{wt}} }
\DeclareMathOperator \Hw  { {\mathrm{hw}} }
\DeclareMathOperator \Supp { {\mathrm{Supp}} }

\renewcommand \max { {\mathrm{max}} }
\newcommand \Spec { {\mathrm{Spec}} }
\begin{document}
%%%%%%%%%%%%%%%%%%%%%%%%%%%%%%%%%%%%%%%%%%%%%%%%%%%%%%%%%%%%%%%%%%%%%%%%%%%
%%%%%%%%%%%%%%%%%%%%%%    Title    %%%%%%%%%%%%%%%%%%%%%%%%%%%%%%%%%%%%%%%%
\title[Multiparameter quantum Schubert cells]
{Spectra and catenarity of multiparameter \\ quantum Schubert cells}
\dedicatory{Dedicated to Kenny Brown and Toby Stafford 
on the occasion of their 60th birthdays}
\author[Milen Yakimov]{Milen Yakimov}
\address{
Department of Mathematics \\
Louisiana State University \\
Baton Rouge, LA 70803
U.S.A.
}
\email{yakimov@math.lsu.edu}
\thanks{The author was supported in part
by NSF grant DMS-1001632.}
\date{}
\begin{abstract}
We study the ring theory of the multiparameter deformations 
of the quantum Schubert cell 
algebras obtained from 2-cocycle twists. This is a large family, which
extends the Artin--Schelter--Tate algebras of twisted quantum matrices.
We classify set theoretically the spectra of all such multiparameter 
quantum Schubert cell algebras, construct each of their prime ideals 
by contracting from explicit normal localizations,
and prove formulas for the 
dimensions of their Goodearl--Letzter strata for
base fields of arbitrary characteristic and 
all deformation parameters that are not roots of unity. 
Furthermore, we prove that the spectra of these algebras 
are normally separated and that all such algebras 
are catenary.
\end{abstract}
\maketitle
%%%%%%%%%%%%%%%%%%%%   Introduction   %%%%%%%%%%%%%%%%%%%%%%%%%%%%%%%%%%%%%%%%
\sectionnew{Introduction}
\lb{intro}
The quantum Schubert cell algebras 
(or quantum nilpotent algebras) form a large
family of iterated Ore extensions defined by  
De Concini, Kac, and Procesi \cite{DKP}, and Lusztig \cite{L}.
They are subalgebras of the negative part of a
quantized universal enveloping algebra $\UU_q(\g)$ 
defined over an arbitrary base field $\KK$ and are 
indexed by the elements of the Weyl group $W$ of $\g$. 
The algebra corresponding to $w \in W$, to be denoted by 
$\UU^w_-$, is a deformation of the universal 
enveloping algebra of the nilpotent Lie algebra 
$\n_- \cap w (\n_+)$, where $\n_\pm$ are 
the niradicals of a pair of opposite Borel subalgebras 
of $\g$. Furthermore, the algebras $\UU^w_-$ are isomorphic 
\cite{Y2} to certain localization of quotients of homogeneous 
coordinate rings of quantum partial flag varieties, similarly 
to the classical case of Schubert cells. They specialize 
\cite{Y3} to the coordinate rings of Schubert cells
equipped with the standard Poisson structure \cite{BGY,GY}.
  
There has been a great interest in the ring theoretic study
of these algebras in an attempt to develop a somewhat general 
theory for iterated Ore extensions. The algebra $\UU^w_-$ 
admits a natural rational action of the torus 
$\Tset^{|\SS(w)|} = (\KK^*)^{\times |\SS(w)|}$
by algebra automorphisms,
where $\SS(w)$ is the support of $w$, see \S 2.4.
The $\Tset^{|\SS(w)|}$-prime spectrum of $\UU^w_-$ was 
described in \cite{MC,Y1}, and the dimensions of the 
corresponding Goodearl--Letzter strata were computed 
in \cite{BCL,Y3}. The Goodearl--Lenagan conjecture \cite{GLen-w}
on polynormality of the torus invariant prime ideals 
of quantum matrices and more generally the algebras 
$\UU^w_-$ was proved in \cite{Y4}, where catenarity and 
normal separation of the spectra of $\UU^w_-$ was established too.    
The algebras $\UU^w_-$ play an important role in various other contexts:
Heckenberger--Schneider classified \cite{HS} the homogeneous coideal 
subalgebras of quantized universal enveloping algebras of Borel subalgebras 
in terms of $\UU^w_-$ and Gei\ss--Leclerc--Schr\"oer proved \cite{GLS}
that $\UU^w_-$ are quantum cluster algebras for simply laced 
$\g$, $\charr \KK=0$, and $q \in \KK$ transcendental over $\Qset$. 

Artin, Schelter, and Tate proved \cite{AST} that multiparameter 
versions of the algebras of quantum matrices and the 
quantum linear groups can be obtained from the 
single parameter case by 2-cocycle twists. 
After this a number of authors investigated 
the effect of such twists on the spectra of graded 
algebras. The spectra of multiparameter quantum groups
was studied by Hodges--Levasseur--Toro \cite{HLT} 
and Costantini--Varagnolo \cite{CV}. The spectra 
of special cases of multiparameter 
quantum Schubert cell algebras were investigated 
by many authors: \cite{GLen-qd,H,O,T}, to name a few. 
Nowlin and Johnson 
\cite{N,JN} proved that certain interesting classes 
of algebras defined in a completely independent way 
turn out to be isomorphic to special twists of 
quantum Schubert cell algebras for affine Kac--Moody algebras.

In this paper we carry out a general study of the ring 
theory of twisted quantum Schubert cell algebras
when $q$ is not a root of unity (without any assumptions 
on the base field $\KK$).
The algebra $\UU^w_-$ is graded by the 
subgroup $\QQ_{\SS(w)}$ of the root lattice of $\g$, 
generated by the simple roots in the support $\SS(w)$ 
of $w \in W$, see \S \ref{2.4}. For 
$\pb \in Z^2(\QQ_{\SS(w)}, \KK^*)$ denote by 
$\UU^w_{-, \pb}$ the algebra obtained by a 
cocycle twist from $\UU^w_-$ using $\pb$. 
Firstly, we give an explicit classification of 
$\Spec \UU^w_{-, \pb}$. This is stated in \thref{spect}.
As in \cite{Y2,Y4} we use 
results of Joseph \cite{J0,J} and Gorelik \cite{Go}. 
Most of those results were formulated in \cite{J,Go} 
for $\KK= k(q)$, $\charr k = 0$. We show that the proofs 
of all such needed results work for an arbitrary $q \in \KK^*$ 
which is not a root of unity and without restrictions 
on the characteristic of $\KK$. At some steps we take shortcuts 
using results of Goodearl and Letzter \cite{GL}, \cite[\S II.6]{BG}. 
Furthermore, \thref{spect} expresses each prime ideal 
of $\UU^w_{-, \pb}$ as a contraction from an 
explicit normal localization of a quotient 
of $\UU^w_{-, \pb}$ by a $\Tset^{|\SS(w)|}$-prime.
These localizations are 
smaller than the ones obtained via the Cauchon 
method of deleting derivations \cite{Ca,MC}.

The spectra $\Spec \UU^w_{-, \pb}$ are partitioned \cite{GL}
into a union of Goodearl--Letzter strata, which are 
isomorphic to the spectra of Laurent polynomial rings. 
We prove an explicit formula for the dimensions of the latter
in \thref{dim}, which works for all 2-cocycles $\pb$, $q \in \KK^*$ not 
a root of unity, and arbitrary base fields $\KK$.
In the special one-parameter case this formula was obtained  
by Bell, Casteels, and Launois \cite{BCL}, and the author \cite{Y3} 
when $\charr \KK =0$ and $q$ is transcendental over $\Qset$.
We give a (very short!) new proof of 
the one-parameter case in \prref{dim2}. 

Furthermore, using results of \cite{Y4} we prove 
in an explicit way that the $\Tset^{|\SS(w)|}$-invariant prime 
ideals of the algebras $\UU^w_{-,\pb}$ are equivariantly 
polynormal. In the special case of multiparameter quantum matrices 
this gives a constructive proof of a conjecture of Brown and Goodearl 
\cite[Conjecture II.10.9]{BG}. Moreover, 
we show that $\Spec \UU^w_{-, \pb}$ are 
normally separated and that all algebras 
$\UU^w_{-, \pb}$ are catenary. This provides a very large 
class of iterated Ore extensions for which Gabber's theorem on 
catenarity of universal enveloping algebras 
of solvable Lie algebras can be extended. 
In a related direction, 
in a forthcoming publication we will prove 
another conjecture of Brown and Goodearl \cite[Conjecture II.10.7]{BG}
that all prime ideals of multiparameter 
quantum groups are completely prime under mild conditions 
of the cocycle twist.
%%%%%%%%%%%%%%%%%%%%%%%%%%%%%%%%%%%%%%%%%%%%%%%%%%%%%%%%%%%%%%%%%%%%%%%%%%%%%%
%%%%%%%%%%%%%%%%%%%%%%%%%%%%%%%%%%%%%%%%%%%%%%%%%%%%%
\sectionnew{Quantum groups, quantum Schubert cells, and their twists}
\lb{qalg}
%%%%%%%%%%%%%%%%
\subsection{}
\label{2.1}
Denote $\Nset= \{0,1, \ldots\}$ and $\Zset_+ = \{ 1, 2, \ldots \}$.
For $m \leq n \in \Zset$ set $[m,n] =\{m, \ldots, n \}$.
Throughout the paper $\KK$ will denote a base field
(of arbitrary characteristic) and $q \in \KK^*$ 
will denote an element which is not 
a root of unity. We fix a simple Lie algebra 
$\g$ of rank $r$ with Cartan matrix $(c_{ij})$. Let $\UU_q(\g)$ be 
the quantized universal enveloping algebra of $\g$ over $\KK$ 
with deformation parameter $q$. Recall \cite{Ja}
that $\UU_q(\g)$ is the $\KK$-algebra with generators 
\[
X^\pm_i, K_i^{\pm 1}, \; \; i \in [1,r]
\]
and relations 
\begin{gather*}
K_i^{-1} K_i = K_i K^{-1}_i = 1, \, K_i K_j = K_j K_i, \,
K_i X^\pm_j K^{-1}_i = q_i^{\pm c_{ij}} X^\pm_j,
\\
X^+_i X^-_j - X^-_j X^+_i = \de_{i,j} \frac{K_i - K^{-1}_i}
{q_i - q^{-1}_i},
\\
\sum_{k=0}^{1-c_{ij}}
(-1)^k
\begin{bmatrix} 
1-c_{ij} \\ k
\end{bmatrix}_{q_{i}}
      (X^\pm_i)^k X^\pm_j (X^\pm_i)^{1-c_{ij}-k} = 0, \, i \neq j,
\end{gather*}
where $q_i = q^{d_i}$ and $\{d_i\}_{i=1}^r$ is the collection of
positive relatively 
prime integers such that $(d_i c_{ij})$ is symmetric. Moreover 
$\UU_q(\g)$ is a Hopf algebra with comultiplication, antipode 
and counit given by
\[
\De(K_i)   = K_i \o K_i, \;
\De(X^+_i) = X^+_i \o 1 + K_i \o X^+_i, \;
\De(X^-_i) = X^-_i \o K_i^{-1} + 1 \o X^-_i
\]
and
\[
S(K_i) = K^{-1}_i, \;
S(X^+_i)= - K^{-1}_i X^+_i, \; 
S(X^-_i)= - X^-_i K_i, \; 
\ep(K_i)=1, \; \ep(X^\pm_i)=0.
\]
The subalgebras of $\UU_q(\g)$ generated by 
$\{X^\pm_1, \ldots, X^\pm_r \}$ will be denoted by 
$\UU_\pm$.
The sets of simple roots, simple coroots, and 
fundamental weights of $\g$ will be denoted by 
$\{\al_i\}_{i=1}^r$, $\{\al_i\spcheck\}_{i=1}^r$,  
and $\{\om_i\}_{i=1}^r$. 
The root and weight lattices of $\g$ will be denoted by 
$\QQ$ and $\PP$. Let
$\QQ^+= \Nset \al_1 + \ldots + \Nset \al_r$, 
$\PP^+ = \Nset \om_1 + \ldots + \Nset \om_r$, and
$\PP^{++} = \Zset_+ \om_1 + \ldots + \Zset_+ \om_r$.
Recall the standard partial order on $\PP$: 
for $\nu_1, \nu_2 \in \PP$ set $\nu_1 \geq \nu_2$ 
if $\nu_2 = \nu_1 - \ga$ for some 
$\ga \in \QQ^+$. Let $\nu_1 > \nu_2$ if $\nu_1 \geq \nu_2$ and
$\nu_1 \neq \nu_2$. If $\la = \sum_{i \in I} n_i \om_i \in \PP^+$, 
$n_i > 0$, $\forall i \in I$ we will say that the support 
of $\la$ is $I$. Denote by $\lcor.,. \rcor$ the symmetric bilinear form on 
$\Span_\Qset \{ \al_i \}_{i=1}^r$ such that 
$\lcor \al_i, \al_j \rcor = d_i c_{ij}$,
$\forall i, j \in [1,r]$. The $q$-weight spaces of a $\UU_q(\g)$-module $V$ 
are defined by 
\[
V_\nu = \{ v \in V \mid K_i v = q^{ \lcor \nu, \al_i \rcor} v, \; 
\forall i \in [1,r] \}, \; \nu \in \PP.
\]
A type one $\UU_q(\g)$-module is a $\UU_q(\g)$-module such that 
$V = \oplus_{\nu \in \PP} V_\nu$.
The category of (left) finite dimensional type one 
$\UU_q(\g)$-modules is semisimple 
(see  \cite[Theorem 5.17]{Ja} and the remark on p. 85 of 
\cite{Ja} for the validity of this for general base fields $\KK$ 
and $q \in \KK^*$ not a root of unity). Furthermore, this category 
is closed under taking tensor products and duals, where the 
latter are defined as left $\UU_q(\g)$-modules using the antipode 
of $\UU_q(\g)$. The irreducible modules in this category  
are parametrized \cite[Theorem 5.10]{Ja} by the set of dominant 
integral weights $\PP^+$. Let $V(\la)$ denote the irreducible 
type one $\UU_q(\g)$-module of highest weight $\la \in \PP^+$. 

Let $W$ and $\B_\g$ be the Weyl and braid groups of $\g$, 
and let $s_1, \ldots, s_r$ and $T_1, \ldots, T_r$ be their 
standard generating sets corresponding to the simple roots 
$\al_1, \ldots, \al_r$. Denote by 
$\ell \colon W \to \Nset$ the standard length function.
The braid group $\B_\g$ acts on 
$\UU_q(\g)$ by algebra automorphisms 
by \cite[eqs. 8.14 (2), (3), (7), and (8)]{Ja} 
and on the finite dimensional type one 
$\UU_q(\g)$-modules by \cite[eq. 8.6 (2)]{Ja}. 
These actions are compatible:
$T_w ( x . v ) = (T_w x) . (T_w v)$ for all 
$w \in W$, $x \in \UU_q(\g)$, $v \in V(\la)$, $\la \in \PP^+$.
%Moreover $T_w(V(\la)_\mu) = V(\la)_{w \mu}$, 
%$\forall w \in W$, $\la \in \PP^+$, $\mu \in \PP$.
%%%%%%%%%%%%%%%%%%%%%%%
\subsection{}
\label{2.2}
If $\KK$ is an algebraically closed field of characteristic 0, 
we will denote by $G$ the connected simply connected algebraic 
group with Lie algebra $\g$. For all base fields $\KK$ and 
deformation parameters $q \in \KK^*$ which are not roots of unity, 
one defines the quantum group $R_q[G]$ as the Hopf subalgebra 
of the restricted dual $(\UU_q(\g))^\circ$ equal to
the span of the matrix coefficients of the modules $V(\la)$, 
$\la \in \PP^+$:
\begin{equation} 
\label{c-notation}
c^\la_{\xi, v} \in (\UU_q(\g))^\circ,\quad 
c^\la_{\xi, v}(x) = \xi ( x v ), \quad v \in V(\la), 
\xi \in V(\la)^*, x \in \UU_q(\g).
\end{equation}
Since we work with arbitrary base fields, $G$ is merely a symbol.
The canonical left and right actions of $\UU_q(\g)$ on 
$\UU_q(\g)^\circ$
\begin{equation}
\label{action}
x \rha c = \sum c_{(2)}(x)c_{(1)}, \; 
c \lha x = \sum c_{(1)}(x)c_{(2)}, \; 
x \in \UU_q(\g), c \in R_q[G]
\end{equation}
preserve $R_q[G]$, where
$\Delta(c) = \sum c_{(1)} \otimes c_{(2)}$.

For each $\la \in \PP^+$ we fix a highest weight 
vector $v_\la$ of $V(\la)$ and denote
\[
c^\la_\xi = c^\la_{\xi, v_\la} \in R_q[G].
\]
The subspace
\[
R^+ = \Span \{ c^\la_\xi \mid \la \in \PP^+, 
\xi \in V(\la)^* \}
\subset R_q[G]
\]
is a subalgebra of $R_q[G]$. We will need the
$R$-matrix commutation relations in $R^+$.
Denote the canonical $\QQ$-grading of $\UU_q(\g)$:
\begin{equation}
\label{Qgrad}
\Wt X^\pm_i = \pm \al_i, \;  \Wt K_i = 0, \quad i  \in [1,r].
\end{equation}
For $\ga \in \QQ^+$, $\ga \neq 0$ set
$m(\ga) = \dim (\UU_+)_\ga= \dim (\UU_-)_{-\ga}$.
Denote by 
$\{u_{\ga, i} \}_{i=1}^{m(\ga)}$ and
$\{u_{-\ga, i} \}_{i=1}^{m(\ga)}$
a pair of dual bases
of $(\UU_+)_\ga$ and $(\UU_-)_{-\ga}$ 
with respect to the Rosso--Tanisaki form
(see \cite[Ch. 6]{Ja} for a discussion 
of the properties of this form 
for arbitrary base fields $\KK$). Then we have:
\begin{equation}
\label{commute}
c_{\xi_1}^{\la_1} c_{\xi_2}^{\la_2} =
q^{ \lcor \la_1, \la_2 \rcor - \lcor \nu_1, \nu_2 \rcor}
c_{\xi_2}^{\la_2} c_{\xi_1}^{\la_1} +
\sum_{\ga \in \QQ^+, \ga \neq 0}
\sum_{i=1}^{m(\ga)}
q^{ \lcor \la_1, \la_2 \rcor - \lcor \nu_1 - \ga , \nu_2 + \ga  \rcor} 
c_{S^{-1}(u_{\ga, i})\xi_2}^{\la_2} 
c_{S^{-1}(u_{-\ga, i}) \xi_1}^{\la_1},
\end{equation}
for all $\la_i \in \PP^+$, $\nu_i \in \PP$, and
$\xi_i \in (V(\la_i)^*)_{\nu_i}$, 
see e.g. \cite[Theorem I.8.15]{BG}.

For $\la \in \PP^+$ and $w \in W$ let
$\xi_{w, \la} \in (V(\la)^*)_{- w\la}$ be such that
$\lcor \xi_{w, \la}, T_w v_\la \rcor =1$. 
(Since $T_w (V(\la)_\la) = V(\la)_{w \la}$, 
$\dim V(\la)_{w \la} = 1$.) Define
\begin{equation}
\label{e}
e^\la_w = c^\la_{\xi_{w,\la}}, \quad \la \in \PP^+, w \in W.
\end{equation}
Then 
\begin{equation}
\label{mult}
e^{\la_1}_w e^{\la_2}_w = e^{\la_1 + \la_2 }_w 
= e^{\la_2}_w e^{\la_1}_w,  \quad \forall \la_1, \la_2 \in \PP^+, 
w \in W,
\end{equation}
see \cite[eq. (2.18)]{Y4}. Denote the multiplicative subsets
$E_w = \{ e_w^\la \mid \la\in \PP^+ \} \subset R^+$.

\ble{Ore} (Joseph, \cite[Lemma 9.1.10]{J})
For all $w \in W$, $E_w$ is an Ore subset 
in $R^+$
\ele
In \cite{J} this result was stated for fields $\KK$ of characteristic 
$0$ (see \cite[\S A.2.9]{J}), but Joseph's proof works for all 
base fields $\KK$, $q \in \KK^*$ not a root 
of unity as we see below. By \eqref{commute} the set $E_{w_0}$
consists of normal elements of $R^+$. Therefore 
$\{ e^n \mid n \in \Nset \}$ is an Ore subset of $R^+$ 
for all $e \in E_{w_0}$, in particular $E_{w_0}$ is 
an Ore subset of $R^+$. Joseph proved 
iteratively that 
\begin{equation}
\label{Orre}
\{e^n \mid n \in \Nset \} \; \; \mbox{are 
Ore subsets of} \; \;  R^+ \; \; 
\mbox{for all} \; \;  e \in E_w, w \in W
\end{equation}
using the following procedure. 
Let $\sig$ be an automorphism of a ring $A$ and 
$\del$ be a locally nilpotent (right skew) $\sig$-derivation 
of $A$ (i.e. for all $a, b \in A$, 
$\del(ab) = \del(a) \sig(b) + a \del(b)$)
such that $\sigma \del \sigma^{-1}$ is a scalar multiple of $\del$.
The degree $\deg_\del a$ of an element $a \in A \backslash \{ 0 \}$ 
is defined as the minimal positive integer $m$ such that 
$\del^{m+1}(a) =0$. Such a skew derivation $\del$ is called 
right regular if for all $\sig$-eigenvectors $a, b \in A$  
of $\del$-degree $m$ and $n$, respectively, and 
$k \in [n, m+n]$
\[
\del^k (a b) = \sum_{i=0}^{m+n-k} s_i ( \del^{k-n+i} a) 
( \del^{n-i} b) 
\]
for some $s_0, \ldots, s_{m+n-k} \in \KK$, $s_0 \neq 0$. 
It is straightforward to check that, 
if $\sig \del \sig^{-1} = q' \del$ for some $q' \in \KK^*$ 
which is not a root of unity, then $\del$ is right regular,
since in that case all coefficients $s$ are products 
of $q'$-binomial coefficients and $\sig$-eigenvalues.  
Joseph's iterative proof of \leref{Ore} relies on the following
fact \cite[Lemma A.2.9]{J}: 

{\em{If $\{ e^n \mid n \in \Nset \}$ 
is an Ore subset of $A$ and $e \in A$ 
is a $\sig$-eigenvector of degree $\deg_\del e = m$, 
then $\{ (\del^m e)^n \mid n \in \Nset \}$ 
is an Ore subset of $A$}}. 

This is applied to
$A = R^+$ and $\sig_i = (\lha K_i^{-1}) $, 
$\del_i = (\lha X_i^-)$, $i \in [1,r]$, 
recall \eqref{action}.
The skew derivations $\del_i$ are regular because 
$\sig_i \del_i \sig_i^{-1} = q^{-2} \del_i$, 
$\forall i \in [1, r]$. 
Let $w s_{i_1} \ldots s_{i_k} = w_0$ where 
$k = \ell(w_0) - \ell(w)$. Denote 
$w_j = w s_{i_1} \ldots s_{i_j}$, $j \in [0, k]$.  
It is easy to show by induction on $j$ that 
$\deg_{\del_{i_j}} e^\la_{w_j} = 
- \lcor w_j \la, \al_{i_j}\spcheck \rcor$, 
$e^\la_{w_{j-1}} = 
t_j \del_{i_{j-1}}^{ - \lcor w_j \la, \al_{i_j}\spcheck \rcor} 
( e^\la_{w_j})$
for some $t_j \in \KK^*$, and that 
$\{ (e^\la_{w_j})^n \mid n \in \Nset \}$ 
is an Ore subset of $R^+$. This implies 
\eqref{Orre} and \leref{Ore}.

Denote the algebras
\[
R^w = R^+ [E_w^{-1}], \quad w \in W
\]
and their invariant subalgebras with respect to 
the action of $\lha K_i$, $i \in [1,r]$,
(cf. \eqref{action}):
\begin{equation}
\label{R0w}
R^w_0 = \{ c^\la_\xi (e^\la_w)^{-1} \mid \la \in \PP^+, \xi \in V(\la)^* \}.
\end{equation}
One does not need to take span in the right hand side of the above formula,
see \cite[eq. (2)]{Go}. 
For $\mu = \la_1 - \la_2 \in \PP$, $\la_1, \la_2 \in \PP^+$
set
\begin{equation}
\label{emu}
e^\mu_w = e^{\la_1}_w (e^{\la_2}_w)^{-1} \in R^+[E_w^{-1}].
\end{equation}
It follows from \eqref{mult} that this does not depend on the choice 
of $\la_1, \la_2$ and that 
$e^{\mu_1}_w e^{\mu_2}_w = e^{\mu_1 + \mu_2}_w$ 
for all $\mu_1, \mu_2 \in \PP$.
For $y \in W$ define
the quantum Schubert cell ideals of $R^+$
\[
Q(y)^\pm = \Span \{ c^\la_\xi \mid \la \in \PP^+, \xi \in V(\la)^*, \,
\xi \perp \UU_\pm T_y v_\la \} 
\]
and the ideals 
\begin{equation}
\label{Qyw}
Q(y)^\pm_w = \{ c^\la_\xi e^{-\la}_w \mid \la \in \PP^+, 
\xi \in V(\la)^*, \,
\xi \perp \UU_\pm T_y v_\la \} = Q(y)^\pm (R^+ [E_w^{-1}]) 
\cap R_0^w 
\end{equation}
of $R_0^w$. Similarly to \eqref{R0w} one does not need 
to take a span in \eqref{Qyw}, see \cite{Go,Y1}.
%%%%%%%%%%%%%%%%%%%%%
\subsection{}
\label{2.3} 
The quantum Schubert cell algebras $\UU^w_\pm$, $w \in W$
were defined by De Concini, Kac, and Procesi 
\cite{DKP}, and Lusztig \cite[\S 40.2]{L} as follows. Let
\begin{equation}
\label{wdecomp}
w = s_{i_1} \ldots s_{i_l}
\end{equation}
be a reduced expression of $w$ (thus $l = \ell(w)$).
Consider the roots
\begin{equation}
\label{beta}
\beta_1 = \al_{i_1}, \beta_2 = s_{i_1} (\al_{i_2}), 
\ldots, \beta_l = s_{i_1} \ldots s_{i_{l-1}} (\al_{i_l})
\end{equation}
and the Lusztig root vectors
\begin{equation}
X^{\pm}_{\beta_1} = X^{\pm}_{i_1}, 
X^{\pm}_{\beta_2} = T_{s_{i_1}} (X^\pm_{i_2}), 
\ldots, X^\pm_{\beta_l} = T_{s_{i_1}} \ldots T_{s_{i_{l-1}}} (X^{\pm}_{i_l}),
\label{rootv}
\end{equation}
see \cite[\S 39.3]{L}. By \cite[Proposition 2.2]{DKP} and 
\cite[Proposition 40.2.1]{L} the subalgebras 
$\UU_\pm^w$ of $\UU_\pm$ generated by 
$X^{\pm}_{\beta_j}$, $j \in [1,l]$ 
do not depend on the choice of a 
reduced decomposition of $w$ and have the PBW bases
\begin{equation}
\label{PBW}
(X^\pm_{\beta_l})^{n_l} \ldots (X^\pm_{\beta_1})^{n_1}, \; \; 
n_1, \ldots, n_l \in \Nset.
\end{equation}

The grading \eqref{Qgrad} 
induces $\QQ$-gradings on the subalgebras $\UU_\pm^w$. 
The corresponding graded components will be denoted by 
$(\UU_\pm^w)_\ga$, $\ga \in \QQ$. The algebra 
$R^+$ is $\PP$-graded by 
\begin{equation}
\label{R+grad}
\Wt c^\la_\xi = \nu, \quad \la \in \PP^+, \nu \in \PP, 
\xi \in (V(\la)^*)_\nu.
\end{equation}
This induces $\PP$-gradings on the algebras $R^w$ and $\QQ$-gradings 
on the algebras $R_0^w$. The latter are given by
\begin{equation}
\label{Qgrad0}
\Wt c^\la_\xi e^{-\la}_w  = \nu + w(\la), \quad
\la \in \PP^+, \nu \in \PP, \xi \in (V(\la)^*)_\nu.
\end{equation}
For a given $\ga \in \QQ^+ \backslash \{ 0 \}$, let
$m_w(\ga) = \dim (\UU^w_+)_\ga= \dim (\UU^w_-)_{-\ga}$
and $\{u_{\ga, i} \}_{i=1}^{m_w(\ga)}$,
$\{u_{-\ga, i} \}_{i=1}^{m_w(\ga)}$ be 
dual bases of $(\UU^w_+)_\ga$ and
$(\UU^w_-)_{-\ga}$ with respect to the Rosso--Tanisaki form,
see \cite[Ch. 6]{Ja}. The quantum $R$ matrix corresponding 
to $w$ is given by
\[
\RR^w = 1 \otimes 1 + \sum_{\ga \in \QQ^+, \ga \neq 0} \sum_{i=1}^{m_w(\ga)} 
u_{\ga, i} \otimes u_{- \ga, i} \in \UU_+ \wh{\otimes} \UU_-,
\]
where $\UU_+ \wh{\otimes} \UU_-$ is 
the completion of $\UU_+ \otimes_\KK \UU_-$ 
with respect to the descending filtration 
\cite[\S 4.1.1]{L}. Recall that there is a unique
algebra antiautomorphism $\tau$ of $\UU_q(\g)$ given by
\begin{equation}
\label{tau}
\tau(X_i^\pm) = X_i^\pm, 
\, 
\tau(K_i) = K_i^{-1}, \; \; 
i = 1, \ldots, r,
\end{equation}
see \cite[Lemma 4.6(b)]{Ja}. It is graded with respect to \eqref{Qgrad}
and satisfies 
\begin{equation}
\label{tau-ident}
\tau (T_w x) = T_{w^{-1}}^{-1} ( \tau (x)), \; \; 
\forall w \in W, x \in \UU_q(\g),
\end{equation}
see \cite[eq. 8.18 (6)]{Ja}. 

\bpr{Q-isom} \cite[Theorem 2.6]{Y4}
Assume that $\KK$ is an arbitrary base field, 
$q \in \KK^*$ is not a root of unity, and $\g$ is a simple Lie 
algebra. For all Weyl group elements $w \in W$ the maps 
\begin{align*}
&\phi_w^+ \colon 
R^w_0 \to \UU^w_-, \quad
\phi_w^+ \big( c^\la_\xi e^{-\la}_w \big) 
= \big( c^\la_{\xi, T_w v_\la} \otimes \id \big) (\tau \otimes \id) \RR^w
\; \; \mbox{and} 
\\
&\phi_w^- \colon 
R^w_0 \to T_w( \UU^{w^{-1} w_0}_+) , \quad
\phi_w^- \big( e^{-\la}_w  c^\la_\xi  \big) 
= \big( \id \otimes c^\la_{\xi, T_w v_\la} \big) \RR^w,
\end{align*}
where $\la \in \PP^+$, $\xi \in V(\la)^*$,
are well defined surjective graded algebra antihomomorphisms 
in the plus case and well defined surjective graded algebra 
homomorphisms in the minus case (with respect to the 
$\QQ$-gradings \eqref{Qgrad} and \eqref{Qgrad0}), 
recall \eqref{c-notation}. Their kernels are $\ker \phi_w^\pm = Q(w)_w^\pm$. 
\epr
The plus case of this proposition was proved in \cite[Theorem 2.6]{Y4}.
The minus case is analogous. In another form the plus case 
was obtained in \cite[Theorem 3.7]{Y1} for a version of $\UU_q(\g)$ 
equpped with the opposite comultiplication, and different 
braid group action and Lustig's root vectors. We note that in the plus and minus cases
the term $e^{-\la}_w$ appears on different sides. This is because of the 
difference in the formulas for the coproducts of $X^\pm_i$. Moreover, the plus case 
requires the use of the antiautomorphism $\tau$, while the minus case 
does not. This is due to the fact that one constructs dual bases of $\UU_\pm$
with respect to the Rosso--Tanisaki form by multiplying Lusztig's 
root vectors in the opposite order,
see e.g. \cite[eqs. 8.30 (1) and (2)]{Ja}.
%%%%%%%%%%%
\subsection{}
\label{2.4}
Let $C$ be an (additive) abelian group and $R$ be a $C$-graded
$\KK$-algebra.  
Following Artin, Schelter, and Tate \cite{AST}, for a 2-cocycle
$\pb \in Z^2(C, \KK^*)$, define a new algebra 
structure on $R$ by twisting the multiplication of $R$ by
\[
b_1 * b_2 = \pb(\al_1, \al_2) b_1 b_2, \quad
\al_i \in C, b_i \in R_{\al_i}, i=1,2.
\]
This algebra will be denoted by $R_\pb$.
Up to an isomorphism $R_\pb$ depends only on the cohomology class 
of $\pb$, see \cite[\S 3]{AST}. Thus we can assume 
that the 2-cocycle $\pb$ is normalized by $\pb(0,0)=1$.
The group of such cocycles will be denoted
by $Z^2(C, \KK^*)_n$. 
The normalization $\pb(0,0)=1$ implies that for all 
$\al \in C$, $\pb(\al,0) = \pb(0, \al) = 1$. Therefore 
the multiplications $R_0 \times R_\al \to R_\al$, 
$R_\al \times R_0 \to R_\al$ remain undeformed for all $\al \in C$
and the multiplicative identity in $R$ is 
an identity of $R_\pb$. Define
\begin{equation}
\label{rb}
\rb(\al, \be) = \pb(\al, \be) \pb(\be, \al)^{-1}, \; \; 
\al, \be \in C.
\end{equation}
By \cite[Proposition 1]{AST}, if $C$ is a free abelian group, then
\begin{equation}
\label{bi-cha}
\rb \colon C \times C \to \KK^* \; \;
\mbox{is a bicharacter}, 
\end{equation}
which is clearly multiplicatively skew symmetric
in the sense that $\rb(\be, \al) = \rb(\al, \be)^{-1}$, 
$\forall \al, \be \in C$.

For $w \in W$ denote the support of $w$:
\[
\SS(w) = \{ i \in [1, r] \mid s_i \leq w \}=
\{ i \in [1, r] \mid s_i \; \; 
{\mbox{appears in a reduced expression of}} \; \; w \}. 
\]
For $I \subseteq [1,r]$ denote 
$\QQ_I = \oplus_{i \in I} \Zset \al_i$,  
$\PP_I = \oplus_{i \in I} \Zset \om_i$, 
$\QQ_I^+ = \QQ_I \cap \QQ^+$,
and $\PP_I^+ = \PP_I \cap \PP^+$. Then
\begin{equation}
\label{Is}
[1,r] \backslash \SS(w) = \{ i \in [1,r] \mid 
w (\om_i) = \om_i \},
\end{equation}
see \cite[Lemma 3.2 and eq. (3.2)]{Y4}. We have:

\ble{grad} \cite[Lemma 3.2 (ii), eq. (2.43)]{Y4} 
The subgroup of $\QQ$ generated by all 
$\ga$ such that $(\UU_\pm^w)_\ga \neq 0$ is equal to $\QQ_{\SS(w)}$.
\ele 

For $\pb \in Z^2(\QQ_{\SS(w)}, \KK^*)_n$ define the 
multiparameter quantum Schubert cell algebra
\[
\UU^w_{\pm, \pb} = (\UU^w_\pm)_\pb. 
\]
For $\ga \in \QQ_{\SS(w)}$ define the character
of $\Tset^{|\SS(w)|}$
\begin{equation}
\label{Tchar}
t \mt t^\ga = \prod_{i \in \SS(w)} t_i^{\lcor \ga, \om_i \rcor}, \quad 
t =(t_i)_{i \in \SS(w)} \in \Tset^{|\SS(w)|}
\end{equation}
and consider the rational $\Tset^{|\SS(w)|}$-action 
on $\UU^w_\pm$ by algebra automorphisms
\begin{equation}
\label{torus-act}
t . x = t^\ga x, \quad 
t \in \Tset^{|\SS(w)|}, 
x \in (\UU^w_\pm)_\ga, \ga \in \QQ_{\SS(w)}.
\end{equation}
Since the twists $\UU^w_{\pm, \pb}$ are graded, 
for all $\pb \in Z^2(\QQ_{\SS(w)}, \KK^*)_n$,
\eqref{torus-act} induces a rational $\Tset^{|\SS(w)|}$-action
on $\UU^w_{\pm, \pb}$ by algebra automorphisms.

For simplicity of the notation define
\begin{equation}
\label{new}
b^\la_{w, \xi} = \phi^+_w(c^\la_\xi e^{-\la}_w), \; 
a^\la_{y, w} = \phi_w^+(e^\la_y e^{-\la}_w ) \in \UU^w_-,
\quad
y \leq w \in W, \la \in \PP^+, \xi \in V(\la)^*
\end{equation}
and denote by $xy=x.y$ 
the multiplication in $\UU^w_{\pm, \pb}$. Eq. \eqref{Is}
implies that $w(\mu) = y(\mu) = \mu$ 
for all $\mu \in \PP_{[1,r]\backslash \SS(w)}$, 
$y \leq w \in W$ and from \eqref{mult} we obtain that 
\begin{equation}
\label{eqq}
a^{\la + \la'}_{y,w} = s a^\la_{y,w}, \quad
\forall \la \in \PP_{\SS(w)}^+,
\la' \in \PP_{[1,r]\backslash \SS(w)}^+, y \leq w \in W,
\end{equation}
for some $s \in \KK^*$.
Recall \prref{Q-isom} and \eqref{new}. Eq. \eqref{commute} 
implies the following commutation relation in 
$\UU^w_{-, \pb}$: 
\begin{align}
\label{commute2}
& \pb(w \la_2 + \nu_2, w \la_1 + \nu_1 )^{-1}
b_{w, \xi_2}^{\la_2} b_{w,\xi_1}^{\la_1} 
=
\pb(w \la_1 + \nu_1, w \la_2 + \nu_2 )^{-1}
q^{ \lcor w \la_1 - \nu_1 , w \la_2 + \nu_2 \rcor }
b_{w, \xi_1}^{\la_1} b_{w, \xi_2}^{\la_2}  +
\\
\nn
&
+ \sum_{\ga \in \QQ^+, \ga \neq 0}
\sum_{i=1}^{m(\ga)}
\pb(w \la_1 + \nu_1 - \ga, w \la_2 + \nu_2 + \ga)^{-1}
q^{ \lcor w \la_1 -\nu_1 +\ga, w \la_2 + \nu_2 + \ga \rcor 
- \lcor \ga, w \la_2 \rcor} \times
\\
&
\nn
\hspace{0.9cm}
\times 
b_{w,S^{-1}(u_{-\ga, i}) \xi_1}^{\la_1}
b_{w, S^{-1}(u_{\ga, i})\xi_2}^{\la_2}
\end{align}
for all $\la_i \in \PP^+$, $\nu_i \in - w \la_i + \QQ_{\SS(w)}$, and
$\xi_i \in (V(\la_i)^*)_{\nu_i}$. The elements 
$u_{\pm \ga, i} \in (\UU_\pm)_{\pm \ga}$ are as 
in \eqref{commute}.
For the values of $\ga \in \QQ^+$ for which 
$w \la_2 + \nu_2 + \ga \notin \QQ(\SS(w))$ or 
$w \la_1 + \nu_1 - \ga \notin \QQ(\SS(w))$, i.e.
$\pb( w \la_2 + \nu_2 + \ga , w \la_1 + \nu_1 - \ga)$ is not defined, 
the term 
$b_{w,S^{-1}(u_{-\ga, i}) \xi_1}^{\la_1}
b_{w, S^{-1}(u_{\ga, i})\xi_2}^{\la_2}$ vanishes
by \prref{Q-isom} and \leref{grad}. The related product in \eqref{commute2}
is set to be equal to $0$. 

The Levendorskii--Soibelman straightening law (see 
for instance \cite[Proposition I.6.10]{BG} and \cite[eq. (2.38)]{Y4}) 
is the following straightening law in $\UU^w_{-, \pb}$: 
\begin{multline}
\label{LS0}
X^-_{\be_j} X^-_{\be_i} - 
\rb(\be_j, \be_i) q^{ - \lcor \be_j, \be_i \rcor }
X^-_{\be_i} X^-_{\be_j}  \\
= \sum_{ {\bf{k}} = (k_{i+1}, \ldots, k_{j-1}) \in \Nset^{\times (j-i-2)} }
s_{\bf{k}} (X^-_{\be_{j-1}})^{k_{j-1}} \ldots (X^-_{\be_{i+1}})^{k_{i+1}},
\; \; s_{\bf{k}} \in \KK,
\end{multline}
for all $i < j$, recall \eqref{rb}. Fix a reduced expression
of $w$ as in  \eqref{wdecomp}. 
For $j \in [0,l]$ denote $w_j= s_{i_1} \ldots s_{i_j}$,  
(thus $w_0=1$ and $w_l = w$). By abuse of notation we will denote 
by the same symbols the restrictions of $\pb$ and $\rb$ to 
$\QQ_{\SS(w_j)} \times \QQ_{\SS(w_j)}$. Because of \eqref{bi-cha},
for all $j \in [1,l]$,
there exists a unique $t_j \in \Tset^{|\SS(w_j)|}$ such that 
\begin{equation}
\label{tj}
t_j^{\beta_i} = 
\rb (\be_j, \be_i) q^{- \lcor \be_j, \be_i \rcor }, \; \; 
\forall i \leq j,
\end{equation}
recall \eqref{Tchar}. The following lemma is a direct consequence
of \eqref{PBW}, \eqref{LS0}, and \eqref{tj}.

\ble{Uw} For all base fields $\KK$, $q \in \KK^*$ not a 
root of unity, Weyl group elements $w \in W$, reduced expressions
\eqref{wdecomp}, and 2-cocycles
$\pb \in Z^2(\QQ_{\SS(w)}, \KK^*)_n$ we have:

(a) The subalgebra of $\UU^w_{-, \pb}$ generated
by $X^-_{\be_1}, \ldots, X^-_{\be_j}$ is isomorphic to 
$\UU^{w_j}_{-, \pb}$, $j \in [1,l]$.

(b) The algebra $\UU^{w_j}_{-, \pb}$ is isomorphic 
to the Ore extension 
$\UU^{w_{j-1}}_{-, \pb}[x, (t_j \cdot), \delta_j]$ 
for some (left) $(t_j \cdot)$-skew derivation $\delta_j$ 
of $\UU^{w_{j-1}}_{-, \pb}$
satisfying $(t_j \cdot) \delta_j = 
q^{- \lcor \be_j, \be_j \rcor} \delta_j (t_j \cdot)$,
$j \in [1,l]$. Here $(t_j \cdot)$ denotes 
the restriction to $\UU^{w_{j-1}}_{-, \pb}$ of the 
automorphism of $\UU^{w_j}_{-, \pb}$ induced from the 
action \eqref{torus-act}, and $\UU^{1}_{-, \pb} = \KK$.

(c) The eigenvalues $t_j . X_{\be_j}^- = 
q^{- \lcor \be_j, \be_j \rcor} X_{\be_j}^-$ are not roots of unity.
\ele

Since the algebras $\UU^w_{-, \pb}$ are iterated Ore
extensions they are noetherian, and the 
Brown--Goodearl theorem \cite[Proposition II.2.9]{BG} implies that each 
$\Tset^{|\SS(w)|}$-prime ideal of $\UU^w_{-, \pb}$ 
with respect to the action \eqref{torus-act} is prime.
%%%%%%%%%%%%%%%%%%%%%%%%%%%%%%%%%%%%%%%%%%%%%
\sectionnew{A description of the spectra of $\UU_{-, \pb}^w$}
\label{spectra}
%%%%%%%%%%%%%%%%%%%%
\subsection{}
\label{3.1}
In this section we describe an explicit stratification of 
$\Spec \UU^w_{-, \pb}$ by the spectra of commutative Laurent 
polynomial rings over $\KK$ and the structure of the poset of 
$\Tset^{|\SS(w)|}$-invariant prime ideals of $\UU^w_{-, \pb}$ 
for the action \eqref{torus-act}.
For $y \in W^{\leq w}$ define the ideals  
\[
I_w(y) = \phi_w^+ ( Q(w)_w^+ + Q(y)_w^-) = 
\phi_w^+ ( Q(y)_w^-)
\]
of $\UU^w_{-, \pb}$ and denote
\[ 
A_{y,w} = \KK^* \{ a_{y,w}^\la \mid \la \in \PP^+ \} =
\KK^* \{ a_{y,w}^\la \mid \la \in \PP^+_{\SS(w)} \},
\]
cf. \eqref{new}--\eqref{eqq}.
By \eqref{mult} and \eqref{commute2}, $A_{y,w}$ is a multiplicative subset 
of $\UU^w_{-, \pb}$ consisting of elements which are normal modulo 
$I_w(y)$:
\begin{equation}
\label{commute3}
a_{y,w}^\la x
=
\rb( (w-y)\la, \ga )
q^{ - \lcor (w+y)\la, \ga \rcor }
x a_{y,w}^\la 
\mod I_w(y),
\end{equation}
$\forall 
\la \in \PP^+,
\ga \in \QQ_{\SS(w)}, 
x \in (\UU_{-, \pb})_\ga$.
The following theorem describes $\Spec \UU^w_{-, \pb}$
set-theoretically and provides some information for its 
Zariski topology.

\bth{spect} Assume that $\KK$ is an arbitrary field, $q \in \KK^*$ is not
a root of unity, $\g$ is a simple Lie algebra, $w$ is an element of 
the corresponding Weyl group $W$, 
and $\pb \in Z^2(\QQ_{\SS(w)}, \KK^*)_n$ is a normalized cocycle.
Then:

(a) The $\Tset^{|\SS(w)|}$-invariant prime ideals of $\UU^w_{-, \pb}$ 
are the ideals $I_w(y)$ for all $y \in W^{\leq w}$. The map 
$y \in W^{\leq w } \mt I_w(y) \in \Tset^{|\SS(w)|}-\Spec \UU^w_{- \pb}$ 
is an isomorphism of posets with respect to the Bruhat order and 
the inclusion order on ideals. In particular, the ideals $I_w(y)$ 
are distinct.

(b) For all $y \in W^{\leq w}$, $I_w(y) \cap A_{y,w} = \emptyset$
and the quotient ring $R_{y,w} = (\UU^w_{-, \pb}/I_w(y))[A_{y,w}^{-1}]$ 
is a $\Tset^{|\SS(w)|}$-simple domain.

(c) For every prime ideal $J$ of $\UU^w_{-, \pb}$ there exists 
a unique $y \in W^{\leq w}$ such that 
\[
J \supseteq I_w(y) \quad
\mbox{and} \quad J \cap A_{y,w} = \emptyset.
\] 
Denote the corresponding subset of $\Spec \UU^w_{-, \pb}$ by 
$\Spec_{I_w(y)} \UU^w_{-, \pb}$, so
\begin{equation}
\label{partition}
\Spec \UU^w_{-, \pb} = \bigsqcup_{y \in W^{\leq w}}
\Spec_{I_w(y)} \UU^w_{-, \pb}.
\end{equation}

(d) The center $Z(R_{y,w})$ is a Laurent polynomial 
ring over $\KK$ and we have a homeomorphism 
\[
\iota_y \colon \Spec Z(R_{y,w}) \to \Spec_{I_w(y)} \UU^w_{-, \pb},
\]
where for $J_0 \in \Spec Z(R_{y,w})$, $\iota_y(J_0)$
is the ideal of $\UU^w_{-, \pb}$ containing 
$I_w(y)$ such that 
$\iota_y(J_0) /I_w(y) = J_0 R_{y,w} \cap (\UU^w_{-, \pb}/I_w(y))$.
The Zariski closures of the strata $\Spec_{I_w(y)} \UU^w_{-, \pb}$ 
are given by
\begin{equation}
\label{closure}
\ol{\Spec_{I_w(y)} \UU^w_{-, \pb}} =
\bigsqcup_{y' \in W^{\leq y}}
\Spec_{I_w(y')} \UU^w_{-, \pb}.
\end{equation}
\eth
Gorelik described in \cite{Go} the torus invariant 
prime ideals of the algebras $R^w_0$ and the resulting 
stratification of $\Spec R^w_0$ as in \thref{spect}
under the condition that $\KK= k(q)$ for a field $k$ 
of characteristic 0. Via the antihomomorphism from 
\prref{Q-isom} this establishes \thref{spect} in the untwisted 
case under those assumptions on $\KK$ and $q$. 
We show that Gorelik's arguments work under the weaker 
assumptions on $\KK$ and $q$, taking shortcuts 
at some steps using results of Goodearl and Letzter \cite{GL, G}.    
The proof of \thref{spect} is given in \S \ref{3.4} and
\S \ref{3.2}--\ref{3.3} contain some 
preparatory results. It follows from 
\thref{spect} (c)-(d) 
that $\Spec_{I_w(y)} \UU^w_{-, \pb}$ are precisely 
the Goodearl--Letzter strata \cite{GL} of $\Spec \UU^w_{-, \pb}$. Namely, 
we have that 
\begin{equation}
\label{GL}
\Spec_{I_w(y)} \UU^w_{-, \pb} = \{ J \in \Spec \UU^w_{-, \pb} \mid
\cap_{t \in \Tset^{|\SS(w)|}} t . J = I_w(y) \}. 
\end{equation}

\bre{restr} For $I \subseteq [1,r]$ 
\[
R^+_I = \Span \{ c^\la_\xi \mid \la \in \PP^+_I, 
\xi \in V(\la)^* \} \subseteq R^+
\]
is a subalgebra of $R^+$ (called the quantum partial flag variety 
associated to $I$), see \cite{Y2}. 
Joseph's argument for the proof of \leref{Ore} shows that 
$E_{\SS(w), w} = \{ e^\la_w \mid \la \in \PP^+_{\SS(w)}\}$
is an Ore subset of $R^+_{\SS(w)}$. Define the subalgebra   
\[
R^w_{\SS(w), 0} = \{ c^\la_\xi (e^\la_w)^{-1} \mid \la \in \PP^+_{\SS(w)}, 
\xi \in V(\la)^* \}
\subset R^+_{\SS(w)} [E_{\SS(w), w}^{-1}]
\]
and its ideals
\begin{equation}
\label{Qrestr}
Q(y)_{\SS(w), w}^\pm =
\{ c^\la_{w,\xi} (e^\la_w)^{-1} 
\mid \la \in \PP^+_{\SS(w)}, \xi \in V(\la)^*, \,
\xi \perp \UU_\pm T_y v_\la \}.
\end{equation}
The argument of the proof of \cite[Theorem 2.6]{Y4} gives 
that the restriction 
\[
\phi^w_+|_{R^w_{\SS(w), 0}} \colon R^w_{\SS(w), 0} \to \UU^w_- 
\]
is a surjective algebra antihomomorphism with kernel 
$Q(w)_{\SS(w), w}^+$, see also \cite[Theorem 3.6]{Y2}.
The proof of \thref{spect} works in the same way 
if the map $\phi^+_w$ is substituted with its restriction 
$\phi^+_w|_{R^w_{\SS(w), 0}}$ and the $\Tset^{|\SS(w)|}$-primes
of $\UU^w_-$ are also given by 
\begin{equation}
\label{Iwy2}
I_w(y) = \phi^w_+( Q(y)_{\SS(w), w}^-), 
\; \; \forall y \in W^{\leq w}.
\end{equation}
\ere
We finish with noting that 
\cite[Theorem 2.3]{GL0} of Goodearl and Letzter, and \leref{Uw}
imply at once:
\bpr{cprime} 
Every prime ideal of $\UU^w_{-, \pb}$ 
is completely prime for all base fields $\KK$, 
$q \in \KK^*$ not a root of unity, Weyl group 
elements $w \in W$, and 2-cocycles 
$\pb \in Z^2(\QQ_{\SS(w)}, \KK^*)_n$
such that the subgroup of $\KK^*$ generated by
\[
\rb(\be_i, \be_j) q^{ \lcor \be_i, \be_j \rcor }
= \left( \rb(\be_j, \be_i) q^{ - \lcor \be_j, \be_i \rcor } \right)^{-1},
\quad 1 \leq i < j \leq \ell(w)
\]
is torsion free, recall \eqref{beta}.
\epr
%%%%%%%%%%
\subsection{}
\label{3.2} We start with a couple of auxiliary 
results for the proof of \thref{spect}.
The following lemma was proved by Joseph
for $\KK = \Cset$, $q \in \Cset^*$ 
not a root of unity \cite[Proposition 10.1.8]{J}
and for $\KK = k(q)$, $\charr k =0$, 
\cite[Proposition 7.3]{J0}.
It is easy to verify that 
Joseph's proof works in the general case. We provide a 
second proof based on \prref{Q-isom}.
\ble{prime1} Let $\KK$ be an arbitrary field, $q \in \KK^*$ 
not a root of unity. The ideals $Q(w)^\pm$ and 
$Q(y)_w^\pm$ of $R^+$ and $R_0^w$ are completely prime
for all $w \in W$, $y \in W^{\leq w}$ in the plus case and
$y \in W^{\geq w}$ in the minus case.  
\ele
\begin{proof}
By \leref{Uw} (b) $\UU^w_\pm$ is a domain, $\forall w \in W$. 
\prref{Q-isom} implies that $Q(w)^\pm_w=\ker \phi^\pm_w$ is 
a completely prime ideal of $R_0^w$. 
It is straightforward to verify that the map
\begin{multline*}
(R_0^w/Q(w)^\pm_w) \# \KK[\PP] \to 
(R^+[E_w^{-1}] / ( Q(w)^\pm R^+[E_w^{-1}]), \\
(c^\la_\xi e^{-\la}_w + Q(w)^\pm_w) \# \mu
\mt c^\la_\xi e^{- \la + \mu}_w + Q(w)^\pm R^+[E_w^{-1}], \; \; 
\mu \in \PP, \la \in \PP^+, \xi \in V(\la)^*
\end{multline*}
is (a well defined) algebra isomorphism, where the smash product 
is defined via the action 
\[
\mu \cdot ( r + Q(w)_w^\pm ) 
= q^{ \mp \lcor w \mu,  \ga   \rcor} ( r + Q(w)_w^\pm ), \quad 
\mu \in \PP, 
r \in (R_0^w)_\ga, \ga \in \QQ, 
\]
cf. \eqref{commute} and \eqref{Qgrad0}.
Therefore $Q(w)^\pm R^+[E_w^{-1}]$ is a completely prime 
ideal of $R^+[E_w^{-1}]$. 
It follows that $Q(w)^\pm =  Q(w)^\pm R^+[E_w^{-1}] \cap R^+$
is a completely prime ideal of $R^+$, $\forall w \in W$.

We have $T_w v_\la \in \UU_\pm T_y v_\la$, for all $\la \in \PP^+$
and $y \in W^{\geq w}$ in the plus case,
$y \in W^{\leq w}$ in the minus case, see \cite[Lemma 4.4.3]{J}. Therefore 
in those cases $Q(y)^\pm \cap E_w = \emptyset$ and thus 
$Q(y)^\pm R^+[E_w^{-1}]$ is a completely prime ideal of $R^+[E_w^{-1}]$.
So $Q(y)_w^\pm = Q(y)^\pm R^+ [E_w^{-1}] \cap R_0^w$ is a 
completely prime ideal of $R_0^w$.     
\end{proof}
\ble{Tprimes} For all base fields $\KK$, $q \in \KK^*$ 
not a root of unity, Weyl group elements 
$y \leq w$, and 2-cocycles
$\pb \in Z^2(\QQ_{\SS(w)}, \KK^*)_n$:

(a) $Q(w)_w^+ + Q(y)_w^-$ is a completely prime ideal of $R^w_0$ 
and 
\begin{equation}
\label{inter}
(Q(w)_w^+ + Q(y)_w^-) \cap \{ e_y^\la e_w^{-\la}
\mid \la \in \PP^+_{\SS(w)} \}  
= \emptyset,
\end{equation}

(b) $I_w(y)= \phi_w^+(Q(w)_w^+ + Q(y)_w^-)$ is a completely prime 
ideal of $\UU^w_{-, \pb}$ and $I_w(y) \cap A_{y, w} = \emptyset$.
\ele

Assume in the setting of \S \ref{2.4} that $C$ is a 
finite rank free abelian group. Then for a graded subspace 
$I$ of $R$, $I$ is a completely prime ideal 
of $R_\pb$ if and only if 
$I$ is a completely prime ideal of $R$, see 
the proof of \cite[Theorem 4.1]{GLen-qd}.
Therefore it is sufficient to prove only 
the untwisted case of \leref{Tprimes} (a). 
Taking into account \prref{Q-isom}, we see that 
part (b) of \leref{Tprimes} follows from the first part of the 
lemma.

Gorelik \cite[Lemma 6.6]{Go} stated the untwisted case of 
\leref{Tprimes} (a) for $\KK = k(q)$, 
$\charr k = 0$. Her proof works under the above more general 
assumptions on $\KK$ and $q$. We review the key steps 
of her proof below to show this. 
For $\mu \in \PP$ denote the automorphism 
\[
\psi_w^\mu(r) = c^{-\mu}_w r c^\mu_w, \quad r \in R^w 
\]
of $R^w$. It is obvious that $R^w_0$ is stable under it.
Recall the skew derivations 
$\del_i = (\lha X_i^-)$ of $R^+$ from \S \ref{2.2}.
By a direct commutation argument one shows that for all $w \in W$, $i \in [1,r]$ 
such that $s_i w < w$, $\la \in \PP^+$, $\mu, \nu \in \PP$, 
$\xi \in (V(\la)^*)_\nu$, and $t$ in the algebraic closure of $\KK$:
\begin{multline}
\label{imp}
\exists k \in \Zset_+ \; \; \mbox{such that} \; \; 
(\psi_w^\mu - t . \id)^k c^\la_\xi = 0 \Rightarrow
\\
\exists k \in \Zset_+ \; \; \mbox{such that} \; \; 
(\psi_w^\mu - t q^{ - \lcor \nu, w \mu \rcor + 
\lcor \nu - m \al_i, s_i w \mu \rcor } \id)^k 
(\del_i^m c^\la_\xi) = 0, 
\end{multline}
where $m \in \Nset$ denotes the $\del_i$-degree 
of $c^\la_\xi$. For $r \in R^w$ set 
\[
\Wt_w(r) = \ga \in \QQ
\; \; \mbox{if} \; \; 
(\psi_w^\mu - q^{ \lcor \ga, \mu \rcor } \id)^k r = 0
\; \; 
\mbox{for some} \; \; 
k \in \Zset_+ \; \; \mbox{and all} 
\; \; \mu \in \PP.
\]
We say that $\xi \in V(\la)^*$
is homogeneous if $\xi \in (V(\la)^*)_\nu$ for some $\nu \in \PP$. 
For a reduced expression $v = s_{i_1} \ldots s_{i_k} \in W$ 
and a homogeneous element $\xi \in V(\la)^*$ denote  
\begin{equation}
\label{tau*}
\del^*_v (c^\la_\xi) : = 
\del_{i_1}^{m_1} \ldots \del_{i_k}^{m_k} (c^\la_\xi),   
\end{equation}
where for $j=k, \ldots, 1$, $m_j$ is 
the $\del_{i_j}$-degree of 
$\del_{i_{j+1}}^{m_{j+1}} \ldots \del_{i_k}^{m_{i_k}} (c^\la_\xi)$,
see \S \ref{2.2}. By induction \eqref{imp} implies
that for all homogeneous $\xi \in V(\la)$, 
if $\Wt_w(c_\xi^\la e_w^{-\la})$ exists, then
\begin{multline}
\label{in}
\Wt_w(c_\xi^\la e_w^{-\la}) 
- w^{-1} \Wt ( c_\xi^\la e_w^{-\la}) 
\\
= \Wt_1 ((\del^*_{w^{-1}}(c^\la_\xi))e_1^{-\la} ) - 
\Wt   ((\del^*_{w^{-1}}(c^\la_\xi)) e_1^{-\la} ) =
- 2 \Wt ((\del^*_{w^{-1}}(c^\la_\xi)) e_1^{-\la} ),
\end{multline}
recall \eqref{mult} and \eqref{Qgrad0}. The last equality follows from 
$e^{-\mu}_1 c^\la_\xi e^\mu_1 = q^{- \lcor \mu, \la+\nu \rcor} c^\la_\xi$, 
$\forall \la \in \PP^+$, $\mu, \nu \in \PP$, $\xi \in (V(\la)^*)_\nu$, 
which is a special case of \eqref{commute}. The automorphisms 
$\psi^\mu_w$, $\mu \in \PP$ of $R^w_0$ commute, because of \eqref{mult} 
and act locally finitely,
since they preserve the grading \eqref{Qgrad0}. Therefore $R^w_0$ is a 
direct sum of the common generalized $\psi^\mu_w$-eigenspaces, $\mu \in \PP$.
Eqs. \eqref{imp} and \eqref{in}, and the fact that 
for each $\ga \in \QQ$ there exists $\la \in \PP^+$ such that 
$(R^w_0)_\ga = \{ c^\la_\xi  e^{-\la}_w \mid \xi \in (V(\la)^*)_{\ga - w \la} \}$ 
imply that
\begin{equation}
\label{Rga}
R_0^w = \bigoplus_{\ga \in - 2 \QQ^+} R_0^w[\ga], 
\end{equation}
where $R_0^w[\ga]$ is the span of all homogeneous elements $r \in R^w_0$ for which 
$\Wt_w(r)$ exists and 
\[
\Wt_w(r) - w^{-1} \Wt(r) = \ga, \; \; \ga \in - 2 \QQ^+.
\]
It follows from \eqref{commute} that
\[
\psi_w^\mu(r) = q^{ \lcor \mu, w^{-1} \ga' \rcor } r \mod Q(w)_w^+, \; \; 
\forall \ga' \in \QQ, r \in (R^w_0)_{\ga'}.
\]
Thus $\oplus_{\ga \in - 2 \QQ^+ \backslash \{ 0 \} } R_0^w[\ga] \subseteq Q(w)_w^+$.
Assume that $R_0^w[0] \cap Q(w)_w^+ \neq 0$. Then there exist $\la \in \PP^+$ 
and a homogeneous element $\xi \in V(\la)^*$, $\xi \neq 0$ such that 
$c_\xi^\la e_w^{-\la} \in R_0^w[0]$ and $\xi( \UU^+ T_w v_\la) =0$.
Because of \eqref{in}, $\Wt ((\del_{w^{-1}}(c^\la_\xi)) e_1^{-\la})=0$,
so $\del_{w^{-1}}(c^\la_\xi) = s e_1^\la$ for some $s \in \KK^*$.
Thus $\xi( (X_{i_1}^-)^{m_1} \ldots (X_{i_l}^-)^{m_l} v_\la) \neq 0$ 
for some $m_1, \ldots, m_l \in \Nset$ where $i_1, \ldots, i_l$ are the 
indices of the reduced decomposition \eqref{wdecomp} of $w$.
Yet 
$(X_{i_1}^-)^{m_1} \ldots (X_{i_l}^-)^{m_l} v_\la \in \UU^+ T_w v_\la$ 
by \cite[Lemma 4.4.3 (v)]{J}, so 
$\xi( (X_{i_1}^-)^{m_1} \ldots (X_{i_l}^-)^{m_l} v_\la) = 0$
because $\xi(\UU_+ T_w v_\la)=0$.
Thus the assumption is not correct and 
\begin{equation}
\label{equal}
Q(w)_w^+ = \oplus_{\ga \in - 2 \QQ^+ \backslash \{ 0 \} } R_0^w[\ga].
\end{equation}
Moreover $R_0^w[0]$ is a subalgebra of $R^w_0$. The ideal
$Q(y)_w^-$ is $\psi^\mu_w$-stable and homogeneous with respect 
to the grading \eqref{Qgrad0}, thus 
$Q(y)_w^- = (Q(y)_w^- \cap R_0^w[0]) \oplus 
(Q(y)_w^- \cap (\oplus_{\ga \in - 2 \QQ^+ \backslash \{ 0 \} } R_0^w[\ga]))$
and
\[
R^w_0/( Q(y)^-_w + Q(w)_w^+) \cong R^w_0[0]/ 
(Q(y)^-_w \cap R^w_0[0]) \cong
(Q(y)^-_w + R^w_0[0])/ Q(y)^-_w \hra
R_0^w/Q(y)^-_w.
\]
Therefore $Q(y)^-_w + Q(w)_w^+$ is a completely prime ideal
of $R^w_0$, since $R^w_0/Q(y)^-_w$ 
is a domain by \leref{prime1}.

Next we go over the key steps of Gorelik's proof of \eqref{inter} 
to show that it works for all base fields $\KK$, $q \in \KK^*$ 
not a root of unity. We have that 
$e^\la_y e^{-\la}_w \in Q(w)^+_w + Q(y)^-_w$ if and only if 
$e^\la_y e^{\la'}_w \in Q(w)^+ + Q(y)^-$ for some $\la' \in \PP^+$. 
(For some of the 
arguments below $\la'$ should be chosen sufficiently large.)
Define $s_i \star w = \max (s_i w, w)$. Let $s_i y >y$ and $y \leq w$. 
Then $s_i y \leq s_i \star w$, \cite[Proposition A.1.7]{J}.
Gorelik proves that 
\begin{equation}
\label{start-ind}
e^\la_y e^{\la'}_w \in Q(w)^+ + Q(y)^- \; \; 
\Rightarrow \; \; 
e^\la_{s_i y} e^{\la'}_{s_i \star w} \in Q(s_i \star w)^+ + Q(s_i y)^-
\end{equation}
as described below.
Since \eqref{inter} is obvious for $y= w = w_0$, 
\eqref{inter} follows from \eqref{start-ind} by induction.

Define the automorphisms $\ol{\sig}_i = (\lha K_i)$
of $R^+$ and the (left skew) derivations   
$\ol{\del}_i = (\lha X_i^+)$, $i \in [1,r]$,
cf. \S \ref{2.2}. 
For a homogeneous element $r \in R^+$ 
with respect to the grading \eqref{R+grad} 
in analogy to \eqref{tau*} define 
$\ol{\del}^*_i(r) = \ol{\del}^n_i (r)$
where $n = \deg_{\ol{\del}_i} r$. 
Since $\KK[X_i^+] \UU_- T_{s_i y} v_\la \subseteq 
\UU_- T_y v_\la$, $\forall \la \in \PP^+$,
\begin{equation}
\label{minus}
\ol{\del}_i^n Q(y)^- \subseteq Q(s_i y)^-, \; \; 
\forall n \in \Nset.  
\end{equation}
If $s_i w < w$, then $\deg_{\ol{\del}_i} e^{\la'}_w = 0$ and 
$e^\la_{s_i y} e^{\la'}_w = 
t \ol{\del}^*_i (e^\la_y e^{\la'}_w ) 
\in Q(w)^+ + Q(s_i y)^-$ for some $t \in \KK^*$, because 
$Q(w)^+$ is $\ol{\del}_i$-invariant. 
This proves \eqref{start-ind} in the case $s_i w < w$.
The case $s_i w > w$ of \eqref{start-ind} is more 
complicated. First one decomposes according to \eqref{Rga}:
\[
e^\la_y e^{-\la}_w = b_0 + \ldots + b_m \; \; 
\mbox{for} \; \; 
b_0 \in Q(y)_w^- \cap R^w_0[0], b_j \in R^w_0[- 2 \ga_j], 
\ga_j \in \QQ^+ \backslash \{0\}, 
j \in [1,m],    
\]
$\ga_j \neq \ga_{j'}$ for $j \neq j'$. Then
\begin{equation}
\label{ww} 
\Wt(b_j) = \Wt( e^\la_y e^{-\la}_w  ) = (w-y) \la, \; \; 
\forall j \in [0,m]. 
\end{equation}
Denote
$n: = \deg_{\ol{\del}^*_i} ( e^\la_y e^{\la'}_w)$.
Since the elements $b_0, \ldots, b_m$ are linearly independent, 
cf. \eqref{Rga}, 
\begin{equation}
\label{Wteqq}
\deg_{\ol{\del}^*_i} ( b_j e^{\la + \la'}_w ) \leq n, 
\; \forall j \in [0,m]
\; \; \mbox{and} \; \; 
e^\la_{s_i y} e^{\la'}_{s_i w}  =
t \ol{\del}^*_i ( e^\la_y e^{\la'}_w ) = 
\sum_{j \in M} t \ol{\del}^*_i ( b_j e^{\la + \la'}_w ), \; t \in \KK^*,
\end{equation}
where $M = \{ j \in [0,m] \mid \deg_{\ol{\del}^*_i} 
( b_j e^{\la + \la'}_w ) = n\}$.
Analogously to \eqref{imp} one shows that for $s_i w >w$, 
if $r \in R^+$ is homogeneous and $\Wt_w (r)$ exists, then 
$\Wt_{s_i w}( \ol{\del}_i r)$ exists and 
\begin{equation}
\label{neweq}
\Wt_w(r) + w^{-1} \Wt(r) = 
\Wt_{s_i w} ( \ol{\del}_i r ) + (s_i w)^{-1} \Wt (\ol{\del}_i r).
\end{equation}
Eqs. \eqref{ww} and \eqref{Wteqq} imply that 
$\Wt(\ol{\del}^*_i ( b_j e^{\la + \la'}_w )) = \Wt( e^\la_{s_i y} e^{\la'}_{s_i w} )$, 
for all $j \in M$. It follows from this, \eqref{neweq}, and \eqref{equal} 
that $\sum_{j \in M \backslash \{ 0 \} } \ol{\del}^*_i ( b_j e^{\la + \la'}_w )
\in Q(s_i w)^+$. Eq. \eqref{minus} 
implies that $\ol{\del}^*_i ( b_0 e^{\la + \la'}_w ) \in Q(s_i y)^-$, 
which completes the proof of \eqref{start-ind} in the case $s_i w >w$.
%%%%%%%%%%%%%%%%%%%%%%
\subsection{}
\label{3.3}
For $J \in \Spec \UU^w_{-, \pb}$ and $\la \in \PP^+$ denote
\[
C_J(\la) = \{ \nu \in \PP \mid \exists \xi \in (V(\la)^*)_\nu, 
b^\la_{w, \xi} \notin J \}.
\]
For all $\la \in \PP^+$, we have $- w \la \in C_J(\la)$ 
since $a^\la_{w,w} = 1 \notin J$.
Denote by $D_J(\la)$ the set of minimal elements of the set $C_J(\la)$.
\ble{min} 
For all base fields $\KK$, $q \in \KK^*$ not a root of unity, 
$w \in W$, $\pb \in Z^2(\QQ_{\SS(w)}, \KK^*)_n$,
and $I \in \Spec \UU^w_{-, \pb}$, there exists a unique 
$y \in W^{\leq w}$ such that $D_J(\la) = \{ - y \la \}$ for all 
$\la \in \PP^+$.
\ele
Gorelik's analogous result \cite[\S 5.2.1]{Go} was formulated 
under the assumption $\KK= k (q)$ for a field $k$, 
$\charr k =0$. This proof works under the more general 
assumptions on $\KK$ and $q$. We sketch this below. 
\\
\noindent
{\em{Proof of \leref{min}.}}
Fix $J \in \Spec \UU^w_{-, \pb}$.
Let $\la \in \PP^+$,
$\nu \in D_J(\la)$, $\xi \in (V(\la)^*)_\nu, 
b^\la_{w, \xi} \notin J$. Applying \eqref{commute2}, we obtain 
that $b^\la_{w, \xi}$ defines a nonzero normal element of 
$\UU^w_{-, \pb}/J$:
\[
b^\la_{w,\xi} x  = 
\rb(\nu + w \la, \ga) q^{ \lcor \nu- w \la, \ga \rcor} 
x b^\la_{w, \xi} \mod J, \; \; 
\forall x \in (\UU^w_{-, \pb})_\ga, \ga \in \QQ_{\SS(w)}.
\] 
Applying this twice and using that $\rb(.,.)$ is 
multiplicatively skew symmetric (see \S \ref{2.4}), 
we obtain that for all
$\la_i \in \PP^+$, $\nu_i \in D_J(\la_i)$,
$\xi_i \in (V(\la_i)^*)_{\nu_i}, 
b^{\la_i}_{w, {\xi_i}} \notin J$, $i=1,2$,
\[
b^{\la_1}_{w,{\xi_1}} b^{\la_2}_{w,{\xi_2}} =
q^{\lcor \nu_1 - w \la_1, \nu_2 + w \la_2 \rcor  
+ \lcor \nu_2 - w \la_2, \nu_1 + w \la_1 \rcor } 
b^{\la_1}_{w,{\xi_1}} b^{\la_2}_{w,{\xi_2}}
\mod J.
\]
Since $J$ is prime and the images of $b^{\la_i}_{w,{\xi_i}}$ in $\UU^w_{-, \pb}/J$
are nonzero normal elements, they are regular. Using the fact that $q \in \KK^*$
is not a root of unity, we obtain
\[
\lcor \nu_1 - w \la_1, \nu_2 + w \la_2 \rcor  
+ \lcor \nu_2 - w \la_2, \nu_1 + w \la_1 \rcor =0, \; \;
\mbox{so} \; \; 
\lcor \nu_1, \nu_2 \rcor = \lcor \la_1, \la_2 \rcor. 
\] 
By \cite[Lemma A.1.17]{J}, if $\la_1, \la_2 \in \PP^{++}$, then 
$\nu_i = - y \la_i$, $i=1,2$, for some $y \in W$. Therefore 
there exists $y \in W$ such that $D_J(\la) = \{ - y \la \}$, 
$\forall \la \in \PP^{++}$. Let $\la \in \PP^{++}$. 
Since $b^\la_{w, \xi} \neq 0$
for some $\xi \in (V(\la)^*)_{- y \la}$ and 
$\ker \phi_w^+ = Q(w)_w^+$ (see \prref{Q-isom}), 
we have $(\UU_+ T_w v_\la)_{y \la} \neq 0$. 
By \cite[Proposition 4.4.5]{J}, $y \leq w$.
Analogously to the proof of 
\cite[Proposition 9.3.8]{J} one shows that 
$D_J(\la) = \{ - y \la \}$, $\forall \la \in \PP^+$
for the same Weyl group element $y$.
\qed
%%%%%%%%%%%%%%%%%%
\subsection{}
\label{3.4} We proceed with the proof of \thref{spect}.
For $y \in W^{\leq w}$ denote
\[
\Spec_y \UU^w_{-, \pb} = \{ J \in \Spec \UU^w_{-, \pb} 
\mid D_J(\la) = \{-y \la \}, \; \forall \la \in \PP^+ \}.
\]
The definition of the sets $D_J(\la)$ 
and \leref{min} imply
\begin{equation}
\label{un}
\Spec \UU^w_{-, \pb} = \bigsqcup_{y \in W^{\leq w}} 
\Spec_y \UU^w_{-, \pb}.
\end{equation}
\noindent
{\em{Proof of \thref{spect}.}}
By \leref{Tprimes} (b), $I_w(y) \in \Spec_y \UU^w_{-, \pb}$, 
$\forall y \in W^{\leq w}$. 
Since $\dim (V(\la))_{y \la}=1$, $\forall \la \in \PP^+, y \in W$, 
each stratum $\Spec_y \UU^w_{-, \pb}$ contains a unique
$\Tset^{|\SS(w)|}$-invariant prime ideal, which would have to be 
precisely the ideal $I_w(y)$. This implies part (a) of the theorem
except the statement for the poset structure of 
$\Tset^{|\SS(w)|}-\Spec \UU^w_{-, \pb}$.
It also implies 
\begin{equation}
\label{GLint}
I_w(y) = \cap_{t \in \Tset^{|\SS(w)|}} t . J, \; \; \forall 
J \in \Spec_y \UU^w_{-, \pb},
\end{equation} 
recall \eqref{torus-act}. Furthermore, we have 
\begin{equation}
\label{int2}
\Spec_y \UU^w_{-, \pb} = \{ J \in \Spec \UU^w_{-, \pb} \mid
J \supseteq I_w(y) \; \; \mbox{and} \; \; 
J \cap A_{y,w} = \emptyset \}.
\end{equation} 
The left inclusion in \eqref{int2} follows from \eqref{GLint} and 
the definition of the sets $D_J(\la)$. 
The equality follows from \eqref{un}.

To complete the proof of part (a) note that $y_1 \leq y_2 \in W^{\leq w}$ 
implies $\UU_- T_{y_1} v_\la \supseteq \UU_- T_{y_2} v_\la$, 
see \cite[Lemma 4.4.3]{J},
therefore $Q(y_1)_w^- \subseteq Q(y_2)_w^-$ and 
$I_w(y_1) = \phi^+_w(Q(y_1)_w^-) \subseteq 
\phi^+_w(Q(y_2)_w^-) = I_w(y_2)$. Now assume that 
$I_w(y_1) \subseteq I_w(y_2)$ for some $y_1, y_2 \in W^{\leq w}$.  
Since $a^\la_{y_2, w} \notin I_w(y_2)$ by \leref{Tprimes}, 
$a^\la_{y_2, w} \notin I_w(y_1)$, $\forall \la \in \PP^+$. 
Therefore $T_{y_2} v_\la \in \UU_- T_{y_1} v_\la$
(because otherwise we would get $c^\la_{y_2} e^{-\la}_w \in Q(y_1)_w^-$ 
and $a^\la_{y_2, w} \in I_w(y_1)$). Applying this to any
$\la \in \PP^{++}$ gives $y_1 \leq y_2$ by \cite[Proposition 4.4.5]{J}.  

Part (c) follows from \eqref{un} and \eqref{int2}. The 
rings $R_{y,w}$ are $\Tset^{|\SS(w)|}$-simple, since 
otherwise a $\Tset^{|\SS(w)|}$-invariant maximal ideal
of $R_{y,w}$ will contract to a $\Tset^{|\SS(w)|}$-invariant
prime ideal in the stratum $\Spec_y \UU^w_{-, \pb}$, which properly contains 
$I_w(y)$. This will contradict with the fact that $I_w(y)$ is 
the only $\Tset^{|\SS(w)|}$-invariant prime ideal in  
$\Spec_y \UU^w_{-, \pb}$. This and \leref{Tprimes} (b) 
prove part (b) of the theorem.

It follows from \leref{Uw} that Goodearl's result \cite[Theorem II.6.4]{BG} 
is applicable to the algebras $\UU^w_{-, \pb}$ and as a consequence of this
all of their $\Tset^{|\SS(w)|}$-prime
ideals are strongly rational. This means that the zero components 
of the centers $R_{y,w}$ with respect to the grading \eqref{Qgrad} 
reduce to scalars: $\left(Z(R_{y,w})\right)_0=\KK$, $\forall y \leq w$. 
Next we apply two results of Goodearl and Letzter: 
\cite[Lemma 6.3 (c)]{GL} implies that 
\begin{equation}
\label{dim01}
\dim (Z(R_{y,w}))_\ga = 0 \; \mbox{or} \; 1, \; \; 
\mbox{for all} \; \; \ga \in \QQ_{\SS(w)}
\end{equation}
and that $Z(R_{y,w})$ are Laurent polynomial rings over $\KK$.   
Furthermore, \cite[Corollary 6.5]{GL} gives that 
contraction and 
extension provide mutually inverse homeomorphisms
between $\Spec Z(R_{y,w})$ and $\Spec R_{y,w}$,
for all $y \leq w$.
Finally, because of \eqref{int2} and general localization facts
$J \mt (J/I_w(y)) R_{y,w}$ is a homeomorphism 
between $\Spec_y \UU^w_{-, \pb}$ 
and $\Spec R_{y,w}$. Eq. \eqref{closure} follows from 
part (a) and \eqref{GLint}. This completes the proof of 
part (d). 
\qed
%%%%%%%%%%%%%%%%%%%%%%%%%%%%%%%%%%%%%%%%%%%%%%%%%%%%%%%
\sectionnew{Dimensions of the Goodearl--Letzter strata 
of $\UU_{-, p}^w$} 
\label{dims}
%%%%%%%%%%%%%%%%%%%%
\subsection{}
\label{4.1} The spectra of the algebras $\UU^w_{-, \pb}$ are
partitioned \eqref{partition} into disjoint unions of the 
Goodearl--Letzter strata $\Spec_{I_w(y)} \UU^w_{-, \pb}$, 
cf. \eqref{GL}. Each of them is
homeomorphic to the spectrum of a Laurent polynomial 
ring over $\KK$, namely $Z(R_{y,w})$, see \thref{spect} (d).
To determine $\Spec \UU^w_{-, \pb}$ set theoretically, 
one needs to solve the problem 
for computing the dimensions of these Laurent polynomial 
rings. In this section we obtain an explicit formula for those 
dimensions. This is done in \thref{dim} and \S \ref{4.2}--\ref{4.3} 
contain some preparatory results.

Let $\mu \in \PP$. If $\mu = \la_1 - \la_2$ for some 
$\la_1, \la_2 \in \PP^+$ with disjoint support, 
define
\begin{equation}
\label{a-mu}
a^\mu_{y,w} = 
a^{\la_1}_{y,w} (a^{\la_2}_{y,w})^{-1} \in R_{y,w}.
\end{equation}
Since $\rb(.,.)$ is a bicharacter, see \eqref{bi-cha}, 
eq. \eqref{commute3} implies that 
\begin{equation}
\label{commute4}
a_{y,w}^\mu x
=
\rb( (w-y)\mu, \ga )
q^{ - \lcor (w+y)\mu, \ga \rcor }
x a_{y,w}^\mu,\quad 
\forall 
\mu \in \PP, \ga \in \QQ_{\SS(w)}, 
x \in (R_{y,w})_\ga.
\end{equation}
Recall \eqref{dim01}. For $y \in W^{\leq w}$ 
denote the sublattice of $\QQ_{\SS(w)}$
\begin{equation}
\label{Zlat}
\ZZ_{y,w} = \{ \ga \in \QQ_{\SS(w)} \mid (Z(R_{y,w}))_\ga \neq 0 \}.
\end{equation}
By \thref{spect} (d)
\begin{equation}
\label{Zeq}
Z(R_{y,w}) \; {\mbox{is a Laurent polynomial ring over $\KK$ 
of dimension equal to}} \; \rank \ZZ_{y,w}.
\end{equation}

The special case for untwisted quantum Schubert cell algebras
of the general dimension result in the below stated \thref{dim} 
was previously obtained by Bell, Casteels, and Launois in \cite{BCL}, 
and the author in \cite{Y3} under the additional condition 
that $\charr \KK =0$ and $q$ is transcendental over $\Qset$.
We give a very simple new derivation of this formula:
 
\bpr{dim2} \cite{BCL,Y3}
For all base fields $\KK$, $q \in \KK^*$ not a 
root of unity, and Weyl group elements $y \leq w$, the 
Goodearl--Letzter strata $\Spec_{I_w(y)} \UU^w_-$ are 
homeomorphic to the spectra Laurent polynomial rings over $\KK$ 
of dimension equal to $\dim \Ker(w+y)$.
\epr
\begin{proof} If 
$\mu \in \Ker_\PP (w+y) = \{ \mu \in \PP \mid (w+y) \mu = 0 \}$,
then by \eqref{commute4}, 
$a^\mu_{y,w} \in (Z(R_{y,w}))_{(w-y)\mu} \backslash \{0\}$. 
Therefore $\ZZ_{y,w} \supseteq (w-y) \Ker_\PP(w+y) = 
2 w \Ker_\PP(w+y)$. So $\rank \ZZ_{y, w} \geq \dim \ker(w+y)$.

Let $z \in (Z(R_{y,w}))_\ga \backslash \{ 0 \}$, $\ga \in \QQ_{\SS(w)}$. 
The central property of $z$ and \eqref{commute4} imply 
\[
a^\mu_{y,w} z = 
q^{ - \lcor (w+y) \mu, \ga \rcor } z a^\mu_{y,w} =
q^{ - \lcor (w+y) \mu, \ga \rcor } a^\mu_{y,w} z
\] 
for all $\mu \in \PP$.
Since $q \in \KK^*$ is not a root of unity, 
$\ga \in ((w+y) \PP)^\perp \cap \QQ$, i.e.
$\ZZ_{y,w} \subseteq ((w+y) \PP)^\perp \cap \QQ$.
Because $\rank ((w+y) \PP)^\perp \cap \QQ = \dim \ker(w+y)$, 
$\rank \ZZ_{y, w} \leq \dim \ker(w+y)$.
Therefore $\rank \ZZ_{y, w} = \dim \ker(w+y)$, and the 
proposition follows from \thref{spect} (d), 
cf. \eqref{Zeq}.
\end{proof}
%%%%%%%%%%%%%%%%%%%%%%%%%%%%%
\subsection{}
\label{4.2} Fix a reduced expression 
\begin{equation}
\label{exp}
s_{i_1} \ldots s_{i_l}
\end{equation}
of $w \in W$ and recall 
\eqref{beta}. We will identify the subexpressions of \eqref{exp} with 
the subsets $D \subseteq [1,l]$. For $D \subset [1,l]$ set 
$s_{i_j}^D = s_{i_j}$, if $j \in D$ and $s_{i_j}^D = 1$ otherwise.
Define $w_{(j)}= s_{i_1} \ldots s_{i_j}$,
$w^D_{(j)}= s_{i_1}^D \ldots s_{i_j}^D$,
$w^D = w^D_{(l)}$ and $\ol{w}^D_{(j)} = s_{i_{j+1}}^D \ldots s_{i_l}^D$,
$j \in [1,l]$. 
A subexpression $D \subseteq [1,l]$ is called a Cauchon 
diagram if for all $j \in [1,l-1]$, $s_{i_j} \ol{w}^D_{(j)} > \ol{w}^D_{(j)}$.
Taking inverses establishes a bijection between the set of those 
and the set of the positive 
subexpressions of Marsh and Rietsch \cite{MR} of
the reverse expression of \eqref{exp}. By \cite[Lemma 3.5]{MR}
for each $y \in W^{\leq w}$, there
exists a unique Cauchon diagram $D \subseteq [1,l]$ such that 
$w^D = y$. Using it, define the lattice
\begin{equation}
\label{Qyw-lat}
\QQ_{y,w} = \sum_{j \in [1,l] \backslash D} \Zset w^D_{(j-1)}(\al_j).
\end{equation}
We will see shortly that this does not depend on the choice of a reduced expression of $w$.
The following lemma provides a second characterization of $\QQ_{y,w}$.

\ble{yw} Let $y \leq w \in W$. For a reduced expression \eqref{exp}
of $w$, let $D \subseteq [1,l]$ be the 
Cauchon diagram such that $w^D =y$.
In the notation \eqref{beta}:
\[
\QQ_{y,w} = \sum_{j \in [1,l] \backslash D} \Zset \be_{i_j}.
\]
\ele

By \cite[Lemma 3.2 (ii)]{Y3}, 
$\QQ_{\SS(w)} = \Zset \be_1 +  \ldots + \Zset \be_l$.
Thus $\QQ_{y,w}$ is a sublattice of $\QQ_{\SS(w)}$.
\leref{yw} follows at once form the fact that 
\[
s_{i_1} \ldots s_{i_j} (\al_{i_{j+1}}) =
s_{i_1} \ldots s_{i_{j-1}} (\al_{i_{j+1}}) \mod 
(\Zset s_{i_1} \ldots s_{i_{j-1}} (\al_{i_j}) ), 
\quad \forall j \in [2,l-1].
\]

For $y \in W^{\leq w}$ denote by 
$\Supp( \UU^w_{-, \pb} / I_w (y) )$ the subgroup of $\QQ_{\SS(w)}$
generated by $\ga \in \QQ_{\SS(w)}^+$ such that 
$(\UU^w_{-, \pb} / I_w (y))_\ga \neq 0$, which is the same as the 
subgroup generated by $\ga \in \QQ_{\SS(w)}$ such that 
$(R_{y,w})_\ga \neq 0$.
Obviously, these supports do not depend on $\pb$.
The next result describes them. (As a side result it also implies 
that $\QQ_{y,w}$ is independent of the choice of a reduced 
expression of $w$.)
\bth{supp} For all base fields $\KK$, $q \in \KK^*$ not a 
root of unity, Weyl group elements $y \leq w$, and 
2-cocycles $\pb \in Z^2(\QQ_{\SS(w)}, \KK^*)_n$
\[
\Supp( \UU^w_{-, \pb} / I_w (y) ) = \QQ_{y,w}.
\]
\eth
\begin{proof} We will first show that 
\begin{equation}
\label{supset}
\Supp( \UU^w_- / I_w (y) ) \supseteq \QQ_{y,w}.
\end{equation}

We will use the fact that, if $y', y'' \in W$ and $i \in [1,r]$ are such that 
$y' s_i > y'$ and $s_i y'' > y''$, then $y' s_i y'' > y' y''$. In terms of the 
reduced expression \eqref{exp}, let $D \subseteq [1,l]$ be the Cauchon diagram
such that $y = w^D$. Denote
$[1,l] \backslash D = \{ j_1 < \ldots < j_{l - |D|} \}$. For $m \in [0, l- |D|]$
define 
\[
D_m = D \sqcup \{ j_1, \ldots j_m\}= D \cup [1,j_m], \; y_m = w^{D_m}.
\]
Then $y_0 = y$, $y_{l- |D|} = w$, and 
$y_m = w_{(j_m)} \ol{w}^D_{(j_m)}= w_{(j_m-1)} s_{i_{j_m}}\ol{w}^D_{(j_m)}$,
$\forall m = 1, \ldots, l- |D|$. Moreover, 
$y_{m-1} = w_{(j_m-1)} \ol{w}^D_{(j_m)}$, $\forall m = 1, \ldots, l- |D|$. 
Since \eqref{exp} is a reduced expression 
$w_{(j_m-1)} s_{i_{j_m}} > w_{(j_m-1)}$. Because $D$ is a Cauchon diagram,
$s_{i_{j_m}}\ol{w}^D_{(j_m)} > \ol{w}^D_{(j_m)}$. By the above mentioned 
fact $y=y_0 < y_1 < \ldots < y_{l - |D|} = w$. By \thref{spect} (a)-(b), 
$a^\la_{y_m, w} \notin I_w(y)$, $\forall m \in [1, l - |D|]$.
Therefore $a^\la_{y_m, w} \in (\UU^w_-/I_w(y))_{(w-y_m)\la}$ 
and  $a^\la_{y_{m-1}, w} \in (\UU^w_-/I_w(y))_{(w-y_{m-1})\la}$.
By an easy computation one obtains that 
$ y_m\la = y_{m-1}\la - \lcor \ol{w}^D_{(j_m)}\la, \al_{i_{j_m}}\spcheck \rcor \be_{j_m}$
in terms of the notation \eqref{beta}. Thus 
\[
\lcor \ol{w}^D_{(j_m)}\la, \al_{i_{j_m}}\spcheck \rcor \be_{j_m} \in 
\Supp ( \UU^w_- / I_w (y) ), \quad \forall \la \in \PP^+, m \in [1, l- |D|].
\]  
Since ${\mathrm{gcd}} \{ \lcor \ol{d}^D_{(j_m)}\la, \al_{i_{j_m}}\spcheck \rcor  \mid 
\la \in \PP^+ \} = 1$,   
\[
\beta_j   \in \Supp( \UU^w_- / I_w (y) ), \quad j \in [1,l] \backslash D.
\]
Now \eqref{supset} follows from \leref{yw}.

M\'eriaux and Cauchon \cite{MC} gave another classification of 
$\Tset^{|\SS(w)|}- \Spec \UU^w _-$, 
associating to each $y \in W^{\leq w}$, 
an ideal $J_w(y) \in \Tset^{|\SS(w)|}- \Spec \UU^w _-$. 
It is based no the Cauchon method of deleting derivations \cite{Ca}, 
which is very different from the one in \cite{Y1}.
It is not known 
yet whether $J_w(y) = I_w(y)$, $\forall y \in W^{\leq w}$. 
M\'eriaux and Cauchon \cite{MC} proved that a certain 
localization of $\UU^w_-/J_w(y)$ by an Ore subset of 
homogeneous elements is isomorphic to a quantum 
torus with generators of weight 
$\{ \be_{i_n} \mid n \in [1,l] 
\backslash D \}$, where $D$ is the Cauchon diagram such that $y = w^D$.
Combining this with \leref{yw} implies that 
$\Supp (\UU^w_- / J_w (y)) = \QQ_{y,w}$, 
$\forall y \in W^{\leq w}$.
Since $\{ I_w(y) \mid y \in W^{\leq w} \} = 
\{ J_w(y) \mid y \in W^{\leq w} \}$ and 
$\Supp (\UU^w_- / I_w (y)) \supseteq \QQ_{y,w}$,
$\forall y \in W^{\leq w}$ we have 
\[
\Supp (\UU^w_- / I_w (y)) = \QQ_{y,w}, 
\quad \forall y \in W^{\leq w}.
\]
\end{proof}
%%%%%%%%%%%%%%%%%%%%%%%%%%%%%%%%%%%%%%
\subsection{}
\label{4.3} 
Assume that $R$ is a $\KK$-algebra graded by a group $C$. We will say that $u \in R$ is 
a {\em{diagonal normal element}} if $u \in R_{\ga_0}$ for some $\ga_0 \in C$ 
and for all $\ga \in C$ there exists $t_\ga \in \KK^*$ 
such that $u r = t_\ga r u$, $\forall r \in R_\ga$. 
If $R$ is a prime ring and $u \neq 0$, 
then $t \colon C \to \KK^*$ is a character.

For all $y \leq w \in W$ there exists $n_{y,w} \in \Zset_+$ such that all 
homomorphisms $\QQ_{y, w} \to \Zset$ have the form 
$\ga \in \QQ_{y,w} \mt \lcor \la, \ga \rcor$ for some 
$\la \in (1/n_{y,w}) \PP$.

\bpr{diag-normal} For all base fields $\KK$, $q \in \KK^*$ not a 
root of unity, Weyl group elements $y \leq w$, 2-cocycles
$\pb \in Z^2(\QQ_{\SS(w)}, \KK^*)_n$, and nonzero 
diagonal normal elements $u \in R_{y,w}$,
there exits $\mu \in (1/2n_{y,w}) \PP$
such that $(w-y) \mu \in \QQ_{\SS(w)}$, 
$u \in (R_{y,w})_{(w-y)\mu}$, 
$\lcor (w+y) \mu, \QQ_{y,w} \rcor \subseteq \Zset$, 
and
\[
u x = \rb( (w-y)\mu, \ga ) q^{ - \lcor (w+y)\mu, \ga \rcor } x u, 
\quad \forall \ga \in \QQ_{y,w}, x \in (R_{y,w})_\ga
\]
in terms of  $n_{y,w} \in \Zset_+$ defined above. 
\epr
\begin{proof} The set of diagonal
 normal elements of an 
algebra graded by an abelian group is invariant under 
twisting. Because of this, one only needs to prove 
the proposition for the algebras $\UU^w_-$, i.e. 
when $\pb$ is trivial.

Let $u \in \UU^w_-/I_w(y)$ be such that 
$u (a^\la_{y,w})^{-1} \in R_{y,w}$ 
is a nonzero diagonal normal element
for some $\la \in \PP^+$.
Then \eqref{commute4} implies that $u$ is 
a nonzero diagonal normal element of $\UU^w_-/I_w(y)$. 
By the same reasoning,
if we establish the proposition for $u$,
then its validity for $u (a^\la_{y,w})^{-1}$ 
will follow. Recall \eqref{exp} and \eqref{beta}.
By Cauchon's method of deleting derivations \cite{Ca}
and the nature of the iterated Ore extension $\UU^w_-$
given by  \leref{Uw}, 
(see \cite{MC}), it follows that there exists a localization 
of $\UU^w_-/I_w(y)$ by an Ore set of homogeneous elements, which is 
graded isomorphic to the $\KK$-quantum torus $\TT$ 
with generators $z_j^{\pm 1}$ (of weights 
$\be_{i_j}$), $j \in D'$ for some subset 
$D' \subseteq [1,l]$ and relations 
$z_j z_{j'} = q^{- \lcor \be_{i_j}, \be_{i_{j'}} \rcor} z_{j'} z_j$, 
$\forall j > j' \in D'$. Since $\TT$ is obtained 
from $\UU^w_-/I_w(y)$ by a localization by homogeneous elements,
$u$ is a diagonal normal element of $\TT$. 
Let $D' = \{ j_1 < \ldots < j_k \}$.
Write $u$ as a sum of monomials 
$t z_{j_1}^{m_1} \ldots z_{j_k}^{m_k}$, 
$t \in \KK^*$, $m_1, \ldots, m_k \in \Zset$.
Let $t z_{j_1}^{m_1} \ldots z_{j_k}^{m_k}$ be one such 
monomial that occurs in $u$.
By \thref{supp} the lattice generated by 
$\be_{i_j}$, $j \in D'$ is equal to $\QQ_{y,w}$.
Since $u$ is a diagonal normal element 
\[
u z_n = q^{ \lcor m_1 \be_{i_{j_1}} + \ldots + m_{n-1} \be_{i_{j_{n-1}}} 
- m_{n+1} \be_{i_{j_{n+1}}} - \ldots - m_k \be_{i_{j_k}}, \be_{i_{j_n}} \rcor }
z_n u, \quad \forall n \in [1,k]
\] 
and 
\[
\be_{i_{j_n}} \in \QQ_{y,w} \mt 
\lcor m_1 \be_{i_{j_1}} + \ldots + m_{n-1} \be_{i_{j_{n-1}}}
- m_{n+1} \be_{i_{j_{n+1}}} - \ldots - m_k \be_{i_{j_k}}, \be_{i_{j_n}} \rcor
\in \Zset, \; \; n \in [1,k] 
\]
defines a group homomorphism. Therefore there exists $\mu_0 \in (1/n_{y,w}) \PP$ 
such that 
\begin{equation}
\label{u-xx}
u x = q^{\lcor \mu_0, \ga \rcor} x u, \; \; 
\forall x \in \TT_\ga, \ga \in \QQ_{y,w}. 
\end{equation}
This implies that \eqref{u-xx} holds for all $x \in (\UU^w_-/I_w(y))_\ga$ 
and thus for all $x \in (R_{y,w})_\ga$. 
Let $u \in (\UU^w_-/I_w(y))_{\ga_0}$, $\ga_0 \in \QQ_{y,w}$.
Using \eqref{commute4} we obtain that 
for all $\la \in \PP$
\[
u a^\la_{y,w} = q^{ \lcor \mu_0, (w-y)\la \rcor } a^\la_{y,w} u = 
q^{ \lcor \mu_0, (w-y)\la \rcor } q^{ - \lcor (w+y)\la, \ga_0 \rcor} u a^\la_{y,w}.
\]
Since $q \in \KK^*$ is not a root of unity 
$\lcor \mu_0, (w-y)\la \rcor - \lcor (w+y)\la, \ga_0 \rcor = 0$, 
i.e.
\[
\lcor \la, w^{-1} (\mu_0 - \ga_0) \rcor = 
\lcor \la, y^{-1} (\mu_0 + \ga_0) \rcor,
\quad \forall \la \in \PP.
\]
Thus $w^{-1} (\mu_0 - \ga_0) = y^{-1} (\mu_0 + \ga_0)$ and 
\[
(w y^{-1} -1) (- \mu_0) = (w y^{-1} +1) \ga_0. 
\]
By standard linear algebra for Cayley transforms there 
exits $\mu \in \Span_\Qset\{\al_i \}_{i=1}^r$ such that 
\[
\ga_0 = (w y^{-1} -1) y \mu = (w -y )\mu \; \; 
\mbox{and} \; \; 
- \mu_0 = (w y^{-1} + 1) y \mu = (w+y )\mu
\]
(see e.g. the proof of \cite[Theorem 3.6]{Y4}).
We have $w \mu = (\ga_0 - \mu_0)/2 \in (1/2 n_{y,w}) \PP$, 
so $\mu \in (1/2 n_{y,w}) \PP$. Furthermore
$(w-y) \mu = \ga_0 \in \QQ_{\SS(w)}$ and 
$\lcor (w+ y) \mu, \QQ_{y,w} \rcor = \lcor
- \mu_0, \QQ_{y,w} \rcor \subseteq \Zset$, 
so $\mu$ satisfies all required properties.
This completes the proof of the proposition.  
\end{proof}
By developing further the arguments of the proof of \prref{diag-normal} 
we will classify all normal elements of the localizations $R_{y,w}$ 
in a forthcoming publication.
%%%%%%%%%%%%%%%%%%%%%%%%%%%%%%%
\subsection{}
\label{4.4} 
Denote the lattice
\begin{equation}
\label{L}
\LL_{y,w, \pb} = \{ (w - y) \mu \mid \mu \in \PP, \;
\rb( (w-y)\mu, \ga )
q^{ - \lcor (w+y)\mu, \ga \rcor } = 1, \; 
\forall \ga \in \QQ_{y,w}
\},
\end{equation}
recall \eqref{Qyw-lat}. One can equivalently use only
$\mu \in \PP_{\SS(w)}$ in \eqref{L} because of \eqref{Is}.

The following theorem provides an explicit formula for the
dimensions of the 
Goodearl--Letzter strata $\Spec_{I_w(y)} \UU^w_{-, \pb}$ 
and completes the set theoretic description 
of the spectra of all multiparameter quantum Schubert cell algebras
$\UU^w_{-, \pb}$.

\bth{dim} For all base fields $\KK$, $q \in \KK^*$ not a 
root of unity, Weyl group elements $y \leq w$, and 2-cocycles
$\pb \in Z^2(\QQ_{\SS(w)}, \KK^*)_n$, the Goodearl--Letzter 
strata $\Spec_{I_w(y)} \UU^w_{-, \pb}$ are homeomorphic to 
the spectra of 
Laurent polynomial rings over $\KK$ of dimension equal to 
the rank of the lattice $\LL^u_{y,w, \pb}$, see \eqref{L}.
\eth
\begin{proof} Recall the definition of $n_{y,w} \in \Zset_+$ from \S \ref{4.3} and 
define the lattice
\begin{align*}
\LL'_{y,w, \pb} = \{ (w - y) \mu \mid &\mu \in (1/ 2n_{y,w})\PP, \;
(w-y) \mu \in \QQ_{\SS(w)}, \; \lcor (w+y) \mu, \QQ_{y,w} \rcor \subseteq \Zset, 
\\
&\rb( (w-y)\mu, \ga )
q^{ - \lcor (w+y)\mu, \ga \rcor } = 1, \; 
\forall \ga \in \QQ_{y,w} \}.
\end{align*}
It follows from \eqref{commute4} that 
$a^\mu_{y,w} \in (\ZZ_{y,w})_{(w-y)\mu} \backslash \{0\}$ for all 
$\mu \in \PP$ that satisfy the equation in \eqref{L}. Therefore 
$\ZZ_{y,w} \supseteq \LL_{y,w, \pb}$. Since every central 
element of a graded ring is a diagonal normal element, 
\prref{diag-normal} implies that 
$\ZZ_{y,w} \subseteq \LL'_{y,w, \pb}$. Therefore
\[
\LL_{y,w, \pb} \subseteq \ZZ_{y,w} \subseteq \LL'_{y,w, \pb}.
\]
One easily shows that the index $[\LL'_{y,w, \pb} : \LL_{y,w, \pb}]$
is finite, starting from $[\PP: (1/2n_{y,w})\PP]< \infty$. Therefore 
\[
\rank \LL_{y,w, \pb} = \rank \ZZ_{y,w} = \rank \LL'_{y,w, \pb}.
\]
The theorem follows from the fact \eqref{Zeq} that $Z(R_{y,w})$ is a 
Laurent polynomial ring of dimension $\rank \ZZ_{y,w}$
and \thref{spect} (d).
\end{proof}
\bre{dim-form} Cauchon \cite{Ca} proved that for a CGL 
extensions $R$ and a torus invariant prime ideal $I$ of $R$, $R/I$ 
admits a localization, which is isomorphic to a quantum torus.
Computing the dimension of the 
corresponding Goodearl--Letzter stratum amounts to 
computing the dimension of its center. Bell and Launois 
applied this to obtain a formula in \cite[Proposition 3.3]{BL} 
for uniparameter CGL extensions, which correspond to 
cases when the Cauchon localization is a special 
kind of quantum torus with exchange relations involving powers of $q$.
For cocycles $\pb$ taking values in the cyclic subgroup 
of $\KK^*$ generated by $q$, this leads to a formula for 
the Goodearl--Letzter strata of $\UU^w_{-, \pb}$
in terms of the dimensions of the kernels of large 
square matrices of size $\ell(w)-\ell(y)$, which is less efficient than 
\thref{dim}.
\ere
%%%%%%%%%%%%%%%%%%%%%%%%%%%%%%%%%%%%%%%%%%%%%%%%%%%%%%%%%%%%%%%%%%%%%%%%%
\sectionnew{Equivariant polynormality, normal separation, and catenarity of 
$\UU^w_{-, \pb}$}
\label{sec5}
%%%%%%%%%%%%%%%
\subsection{}
\label{5.1}
In this section, using results of our previous paper \cite{Y5}, 
we prove that all $\Tset^{|\SS(w)|}$-prime ideals of $\UU^w_{-, \pb}$ 
are equivariantly polynormal,
that the spectra $\Spec \UU^w_{-, \pb}$ are normally separated,
and that all algebras $\UU^w_{-, \pb}$ are catenary.

We start with a brief review of equivariant polynormality. Assume that
a ring $R$ is equipped with an action of a group $\Ga$
by algebra automorphism.
We say that an element $u \in R$ is $\Ga$-normal if it is a 
$\Ga$-eigenvector and if there exists $g \in \Ga$ such that 
$u r = (g . r) u$ for all $r \in R$. Sometimes $\Ga$-normality is 
defined requiring only the second condition, 
see \cite{G}, but this is not sufficient to extend the 
definition to $\Ga$-polynormality for $\Ga$-stable ideals of $R$.
Given a $\Ga$-stable ideal $I$ of $R$, we say that 
a sequence $u_1, \ldots, u_N \in R$ is a $\Ga$-polynormal 
generating sequence of $I$ if $\{u_1, \ldots, u_N\}$ generates $I$ 
and for all $i = 1, \ldots, N$ the element $u_i$ is a $\Ga$-normal 
element of $R$ modulo the ideal generated by $u_1, \ldots, u_{i-1}$. 
Note that it follows from the the conditions on the elements 
$u_1, \ldots, u_{i-1}$ that the two-sided ideal of $R$ generated by them 
is $\Ga$-stable and equals $R u_1 + \ldots + R u_i$. 
In particular, $R$ is generated by $u_1, \ldots, u_N$ as a 
one-sided (left or right) ideal. 

For all $\mu \in \PP_{\SS(w)}$, $\vartheta \in \QQ_{\SS(w)}$ there exists a unique 
$t_{\mu, \vartheta} \in \Tset^{|\SS(w)|}$ such that
\[
\left( t_{\mu, \vartheta} \right)^\ga =
\rb(\vartheta, \ga )
q^{ \lcor \vartheta- 2 w \mu,  \ga \rcor }, \; \; \forall \ga \in \QQ_{\SS(w)},
\]
cf. \eqref{Tchar} and \eqref{commute2}. It is clear that 
the map $(\mu, \vartheta) \in \PP_{\SS(w)} \times \QQ_{\SS(w)} \mt t_{\mu, \vartheta} 
\in \Tset^{|\SS(w)|}$ is a group homomorphism. 
Thus its image
\begin{equation}
\label{PQpb}
\PP\QQ_{w, \pb} = \{ t_{\mu, \vartheta} \mid 
\mu \in \PP_{\SS(w)}, \vartheta \in \QQ_{\SS(w)} \}
\end{equation}
is a subgroup of $\Tset^{|\SS(w)|}$,
which is a quotient of $\PP_{\SS(w)} \times \QQ_{\SS(w)}$.

For every $\la \in \PP^+_{\SS(w)}$
and $y \leq w \in W$ fix a finite subset $\VV_{y,w}(\la)^* \subset V(\la)^*$ 
consisting of homogeneous elements, such that the restriction map 
\[
\xi \in V(\la)^* \mt \xi|_{\UU_+ T_w v_\la} \in (\UU_+ T_w v_\la)^*
\]
sends bijectively $\VV_{y,w}(\la)^*$ to a basis of 
$(\UU_- T_y v_\la \cap \UU_+ T_w v_\la)^\perp 
\subset ( \UU_+ T_w v_\la)^*$. Given $\Om \subset \PP^+_{\SS(w)}$ denote 
\[
\VV_{y,w}(\Om) = \bigsqcup_{\la \in \Om} \VV_{y,w}(\la).
\]
Define the maps
\[
\Hw \colon \VV_{y,w}(\Om) \to \Om, 
\Wt \colon \VV_{y,w}(\om) \to \PP, \quad
\Hw (\xi) = \la, \Wt (\xi) = \nu, \; \; 
\mbox{if} \; \; \xi \in (V(\la)^*)_\nu. 
\] 
Denote the partial ordering on $\VV_{y,w}(\Om)$:
\begin{equation}
\label{po}
\xi \preceq \xi', \; \; \mbox{if} \; \; \Hw(\xi) = \Hw(\xi') \; \; 
\mbox{and} \; \; \Wt(\xi) \leq \Wt(\xi').
\end{equation}

\bth{polynorm} Assume that $\KK$ is an arbitrary base field, 
$q \in \KK^*$ is not a root of unity, $\g$ is an arbitrary 
simple Lie algebra, $w \in W$, 
and $\pb \in Z^2(\QQ_{\SS(w)}, \KK^*)_n$ is a normalized 2-cocycle.

(a) Since $\UU^w_{-, \pb}$ is noetherian, for each $y \in W^{\leq w}$ there 
exists a finite subset $\Om$ of $\PP^+_{\SS(w)}$ such that 
the ideal $I_w(y)$ is generated by $b^\la_{w, \xi}$ for $\la \in \Om$, 
$\xi \in V(\la)^*$, recall \eqref{Iwy2}. Fix 
any linear ordering $\xi_1 < \ldots < \xi_N$ on $\VV_{y,w}(\Om)$,
which is a refinement of the partial ordering \eqref{po}. Then 
\[
b^{\Hw(\xi_1)}_{w, \xi_1}, \ldots, b^{\Hw(\xi_N)}_{w, \xi_N}
\]
is a $\Tset^{|\SS(w)|}$-polynormal generating sequence of the ideal $I_w(y)$. 
More precisely 
\[
b^{\Hw(\xi_i)}_{w, \xi_i} x =
(t_{\Hw(\xi_i), \Wt(\xi_i) + w(\Hw(\xi_i))} . x) x b^{\Hw(\xi)}_{w, \xi_i}
\mod 
\left( \UU^w_{-, \pb} b^{\Hw(\xi_1)}_{w, \xi_1} + \ldots +
\UU^w_{-, \pb} b^{\Hw(\xi_{i-1})}_{w, \xi_{i-1}} 
\right)
\]
for all $x \in \UU^w_{-, \pb}$, $i=1, \ldots, N$. 
All automorphism related to the above normal elements 
come from the subgroup $\PP\QQ_{w, \pb}$ 
of $\Tset^{|\SS(w)|}$, cf. \eqref{PQpb}.

(b) If the base field $\KK$ has characteristic $0$ and 
$q$ is transcendental over $\Qset$, then the 
conclusion of part (a) is valid for $\Om = \{\om_i \mid i \in \SS(w) \}$.
\eth
\begin{proof} The untwisted cases of parts (a) and (b) of
\thref{polynorm} were proved in \cite[Theorem 3.6]{Y5} and
\cite[Theorem 3.4]{Y5}, respectively. The theorem follows from this, 
since all elements of the polynormal generating sets 
are homogeneous with respect to the grading \eqref{Qgrad}.
\end{proof} 
A special case of \thref{polynorm} gives a constructive proof of 
a conjecture of Brown and Goodearl, \cite[Conjecture II.10.9]{BG},
that the torus invariant prime ideals of the multiparameter 
algebras of quantum matrices have polynormal generating sets
consisting of quantum minors. 
(This is an extension of the Goodearl--Lenagan conjecture \cite{GLen-w}
on polynormality in single parameter quantum matrix algebras,
which we proved in \cite{Y5}.) Artin, Schelter, and Tate showed 
\cite[p. 889]{AST} (see also \cite[Lemma 3.6]{GLen-qd}) that
the multiparameter algebras of quantum matrices 
${\mathcal{O}}_{\la, \pb}(M_{m,n}(\KK))$ are obtained by twists 
from the single parameter ones if $\la$ is a square root of $q$ in 
$\KK$ (provided that such exists). In \cite[Proposition 2.1.1]{MC} 
(and \cite[Lemma 4.1]{Y5}) the single parameter algebras of quantum 
matrices were realized as special cases of the algebras $\UU^w_+$ 
(and $\UU^w_-$) in a way that matches the corresponding gradings 
by free abelian groups. 
Therefore the isomorphism of \cite[Lemma 4.1]{Y5}
realizes ${\mathcal{O}}_{\la, \pb}(M_{m,n}(\KK))$ as a special case 
of the algebras $\UU^w_{-, \pb}$. In this case 
by \cite[Lemma 4.3]{Y5},
if $\charr \KK =0$ and $q$ is transcendental over $\Qset$, 
the generating sets from \thref{polynorm} (b) consist 
of quantum minors, which are explicitly listed 
in \cite[Theorem 4.4]{Y5}. More generally, for all simple 
Lie algebras $\g$ and Weyl group elements $w$,
if $\charr \KK =0$ and $q$ is transcendental over $\Qset$,
then in \thref{polynorm} one can choose $\Om=\{\om_1, \ldots, \om_r\}$
by \cite[Theoerm 3.4]{Y5}. The last two facts rely on 
a theorem of Joseph \cite[Th\'eor\`eme 3]{J2} and need the stronger
assumption on $\KK$ and $q$ because of a specialization argument. 
%%%%%%%%%%%%%%%%%%
\subsection{}
\label{5.2}
Finally we establish normal separation of $\Spec \UU^w_{-, \pb}$ 
and catenarity of $\UU^w_{-, \pb}$.
 
\bth{nsepar} For all base fields $\KK$, $q \in \KK^*$ not a root of unity, 
simple Lie algebras $\g$, $w \in W$, and 2-cocycles
$\pb \in Z^2(\QQ_{\SS(w)}, \KK^*)_n$,
the spectra of the algebras $\UU^w_{-, \pb}$ are normally separated.
\eth
\begin{proof} Using Goodearl's result \cite[Corollary 4.6]{G} it suffices 
to show that $\UU^w_{-, \pb}$ has graded normal separation. This can 
be proved in two different ways. Firstly, assume that 
$I_1 \subsetneq I_2$ are two graded (i.e. $\Tset^{|\SS(w)|}$-invariant) 
prime ideals of $\UU^w_{-, \pb}$. If $x_1, \ldots, x_n$ is a 
$\Tset^{|\SS(w)|}$-polynormal generating sequence of $I_2$
as in \thref{polynorm} (a), 
then the first $x_i$ which does not belong to $I_1$ produces a graded 
normal separating element for the pair of ideals $I_1, I_2$. Secondly, 
by \thref{spect} (a) for every pair $I_1 \subsetneq I_2$ of graded prime
ideals of $\UU^w_{-, \pb}$, there exist $y_1 < y_2 \in W^{\leq w}$ 
such that $I_1 = I_w(y_1)$, $I_2 = I_w(y_2)$. Let $\la \in \PP^{++}$. 
Then $y_1 \la > y_2 \la$, so $a_{y_1, w}^\la \in I_w(y_2)$.
\thref{spect} (c) implies that $a_{y_1, w}^\la \notin I_w(y_1)$ and 
eq. \eqref{commute3} gives that $a_{y_1, w}^\la$ is normal modulo $I_w(y_1)$. 
\end{proof}

Recall that a ring $R$ is Auslander--Gorenstein if 
the injective dimension of $R$ (as both right and left 
$R$-module) is finite, and for all integers $0 \leq i < j$ 
and finitely generated (right or left) $R$-modules $M$, 
we have $\Ext^i_R(N, R) = 0$ 
for all $R$-submodules $N$ of $\Ext^j_R(M,R)$. A ring $R$ is 
Auslander regular if, in addition, the global dimension 
of $R$ is finite. The grade of a finitely generated
$R$-module $M$ is given by
\[
j(M) = \inf \{i \geq 0 \mid \Ext^i_R(M,R) \neq 0 \}.
\]
An algebra $R$ is Cohen--Macauley if 
\[
j(M)+ \GKdim M = \GKdim R 
\]
for all finitely generated $R$-modules $M$. We will need 
the following two results:

\bth{G-L} (Goodearl--Lenagan, \cite{GLen-c}) Assume that $A$ is an 
affine, noetherian, Auslander--Gorenstein and Cohen--Macaulay
algebra over a field, with finite Gelfand--Kirillov dimension. 
If $\Spec A$ is normally separated, then $A$ is catenary. 
If, in addition, $A$ is a prime ring, then Tauvel's height 
formula holds.
\eth
Tauvel's height formula holds for a ring $R$ if for all prime 
ideals $J$ of $R$, the height of $J$ is equal to 
\[
\GKdim R - \GKdim(R/J).
\]

\bth{LS} (Ekstr\"om, Levasseur--Stafford, \cite{E, LS}) Assume $R$ is a 
noetherian, Auslander regular ring. 
Let $S= R[x; \sigma, \delta]$ be an Ore 
extension of $R$. Then: 

(a) $S$ is Auslander regular.

(b) If $R = \oplus_{k \geq 0} R_k$ is 
a connected graded Cohen--Macauley 
$\KK$-algebra over a field $\KK$ such 
that $\sigma (R_k) \subseteq R_k$ 
for all $k \geq 0$, then $S$ is Cohen--Macauley. 
\eth 

The next theorem proves that all algebras $\UU^w_{-, \pb}$ are catenary.

\bth{catenary} For all base fields $\KK$, $q \in \KK^*$ 
not a root of unity, simple Lie algebras $\g$, 
$w \in W$, and 
and 2-cocycles $\pb \in Z^2(\QQ_{\SS(w)}, \KK^*)_n$,  
the algebras  $\UU^w_{-, \pb}$ are Auslander regular, 
Cohen--Macauley, catenary,
and Tauvel's height formula holds for them.
For all $y \in W^{\leq w}$ the height 
of prime ideal $I_w(y)$ is equal to $\ell(y)$ and
\[
\GKdim (\UU^w_{-, \pb}/I_w(y)) = \ell(w) - \ell(y).
\]
\eth
\begin{proof} \leref{Uw} (a)-(b) and 
\thref{LS} (a) imply that $\UU^w_{-, \pb}$ 
is Auslander regular. Given $\la \in \PP^{++}$, we 
can specialize the $-\QQ^+_{\SS(w)}$-grading
of $\UU^w_{-, \pb}$ from \eqref{Qgrad} to an 
$\Nset$-grading by
\[
(\UU^w_{-, \pb})_n = \{ x \in (\UU^w_{-, \pb })_{-\ga} \mid
\ga \in \QQ^+_{\SS(w)}, \lcor \la, \ga \rcor = n \}, \; \; 
n \in \Nset.
\]
Obviously the $\Tset^{|\SS(w)|}$-action \eqref{torus-act} preserves 
each graded component $(\UU^w_{-, \pb})_n$ and 
$\UU^w_{-, \pb}$ is connected. \leref{Uw} (a)-(b) 
and \thref{LS} (b) imply that $\UU^w_{-, \pb}$ is
Cohen--Macauley.
It also follows from \leref{Uw} (a)-(b) that the algebras 
$\UU^w_{-, \pb}$ are affine, noetherian domains 
with $\GKdim \UU^w_{-, \pb} = \ell(w)$. \thref{G-L} implies 
that all algebras $\UU^w_{-, \pb}$ are catenary and 
satisfy Tauvel's height formula.

Since the Gelfand--Kirillov dimension of a graded algebra does not 
change under twisting, the last statement of the theorem 
follows from \cite[Theorem 5.8]{Y5}.
\end{proof}
Some special cases of the untwisted case of Theorems \ref{tpolynorm}, \ref{tnsepar}, 
and \ref{tcatenary} were established by Caldero \cite{Cal}, Cauchon \cite{Ca2}, 
Goodearl--Lenagan \cite{GLen-c}, and Malliavin \cite{M}. The general 
untwisted case of Theorems \ref{tpolynorm}, \ref{tnsepar}, and \ref{tcatenary}
was proved by the author in \cite{Y5}.
%%%%%%%%%%%%%%%%%%%%%% References %%%%%%%%%%%%%%%%%%%%%%%%%%%%%%%%%%%%%%%

%%%%%%%%%%%%%%%%%%%%%%%%%%%%%%%%%%%%%%%%%%%%%%%%%%%%%%%%%%%%%%%%%%%%%%%%%%%%%%%
%%%%%%%%%%%%%%%%%%%%%%%%%%%%%%%%%%%%%%%%%%%%%%%%%%%%%%%%%%%%%%%%%%%%%%%%%%%%%%
\end{document}